\documentclass[11pt,reqno]{amsart}
\usepackage{amsmath,amssymb,amsthm,amscd,latexsym}
\usepackage{graphicx}
\usepackage{color}
\usepackage{tikz}
\usepackage{multirow}
\usepackage{comment}
\usepackage[top=30mm,bottom=30mm,left=30mm,right=30mm]{geometry}

\usepackage{enumitem}

\usetikzlibrary{patterns}
\usetikzlibrary{intersections, calc, decorations.markings}

\usetikzlibrary{decorations.pathmorphing, decorations.pathreplacing, shapes.geometric}

\tikzset{myQV/.style={circle, fill, inner sep=1.5pt, outer sep = 0.5pt}} 
\tikzset{myCN/.style={circle, inner sep = 0, outer sep = 0.5pt}}  %%% circle node
\tikzset{myMQV/.style={star, star points=5, star point ratio=2.25, fill=black, inner sep=1.5pt, outer sep=0.5pt}} %%%% quiver star vertex

\tikzset{>=stealth}

\usepackage[colorlinks=true, linkcolor=blue, citecolor=red]{hyperref}
\usepackage{bookmark}

\numberwithin{equation}{section}
\theoremstyle{definition}
\newtheorem{example}{Example}[section]
\newtheorem{definition}[example]{Definition}

\theoremstyle{plain}
\newtheorem{lemma}[example]{Lemma}
\newtheorem{theorem}[example]{Theorem}
\newtheorem{proposition}[example]{Proposition}
\newtheorem{corollary}[example]{Corollary}
\newtheorem{conjecture}[example]{Conjecture}

\newtheorem{assumption}{Assumption}

\newtheorem{maintheorem}{Theorem}

\newlength\squareheight
\newcommand{\relmiddle}[1]{\mathrel{}\middle#1\mathrel{}}

\renewcommand{\int}{\operatorname{int}}

\DeclareMathOperator{\gseed}{{\rm g-CL}}

\newcommand{\D}{\mathbf{D}}
\newcommand{\bra}[1]{\left[#1\right]}

\title[Finiteness and completeness of $g$-fans]
{Finiteness and completeness of $g$-fans\\ for real-weighted quivers}
\author{Ryota Akagi}
\address{Graduate School of Mathematics, Nagoya University, Chikusa-ku, Nagoya, 464-0813, Japan}
\email{ryota.akagi.e6@math.nagoya-u.ac.jp}
\author{Toshiya Yurikusa}
\address{Department of Mathematics, Osaka Metropolitan University, Osaka 558-8585, Japan}
\email{yurikusa@omu.ac.jp}

\begin{document}

\begin{abstract}
	We study finiteness and completeness of $g$-fans for quivers with positive real weights. We characterize finite Coxeter type by the condition that every weight in the mutation class has the form $2\cos(\pi/m)$ for some integer $m\geq3$. Under total sign-coherence, we then show that a quiver with finitely many $g$-clusters is $2$-finite, that is, every weight in its mutation class is less than $2$. This also gives a new proof that $g$-finiteness implies finite Coxeter type. Assuming in addition that, for every initial quiver in the mutation class, each $c$-vector on a coordinate axis is a signed standard basis vector, we prove that the $g$-fan is complete if and only if there are finitely many $g$-clusters. The main geometric result states that the coordinate projection of every $g$-cone is contained in a $g$-cone of the corresponding full subquiver. Consequently, completeness and density pass to full subquivers.
\end{abstract}

\maketitle

\section{Introduction}\label{sec:intro}
	Cluster algebras were introduced by Fomin and Zelevinsky \cite{FZ02} to provide a combinatorial framework for studying Lusztig's canonical bases in quantum groups and total positivity. Since their inception, cluster structures have been uncovered across a wide variety of mathematical areas, including representation theory, algebraic geometry, and integrable systems.
	
	The heart of cluster algebra theory is a cluster pattern, an inductive system producing distinguished generators termed cluster variables. A cluster pattern consists of seeds, which are transformed via an elementary operation called mutation. While cluster variables were originally the primary objects of study, Fomin and Zelevinsky \cite{FZ07} introduced $c$-vectors and $g$-vectors to capture the discrete and combinatorial information underlying seed mutations. Subsequent developments have revealed that these vectors govern the structural behavior of cluster patterns and share deep analogies with root systems.
	
	For instance, every $c$-vector is sign-coherent; namely, all of its nonzero entries have the same sign \cite{GHKK18}. In addition, the $g$-vectors in each seed span a simplicial cone, and these cones together with their faces form the $g$-fan, providing a geometric realization of the cluster complex \cite{GHKK18,Re14}. In the skew-symmetric case, $c$-vectors naturally appear as normal vectors pointing inward across the facets of maximal cones \cite{NZ12}.

  In particular, $g$-fans themselves exhibit rich combinatorial and geometric structures. In ordinary (integer) cluster algebras, they are seen as analogues of Weyl chamber structures. For example, for certain classes, $g$-fans are obtained by taking unions of some chambers of the corresponding Weyl chamber structure \cite{RS16, RS18}. Furthermore, $g$-fans provide the following characterization of finite type: the cluster algebra is of finite type if and only if the $g$-fan is complete, i.e., it covers the whole space $\mathbb{R}^n$ \cite{N26,Y25}.
	
	Although traditionally defined for integer matrices, the mutation formulas for exchange matrices, $c$-vectors, and $g$-vectors remain well-defined for matrices with real entries. In recent years, several authors have extended the cluster algebraic framework to this real-valued setting. Lampe \cite{Lam18} considered this generalization in the context of noncrystallographic root systems. Subsequently, Duffield and Tumarkin \cite{DT24, DT25} realized $c$-vectors and $g$-vectors for noncrystallographic types via a weighted folding technique. From the perspective of geometric realizations of quiver mutations, Felikson and Tumarkin \cite{FT19} constructed geometric models for mutations of rank $3$ quivers with real weights. Using such realizations, Felikson and Lampe \cite{FL23} studied exchange graphs for non-integer quivers of finite type in arbitrary rank. Furthermore, Felikson and Tumarkin \cite{FT23} gave a complete classification of mutation-finite quivers with real weights, including many examples that do not occur in the classical theory.
	
	Akagi and Chen further developed this framework in \cite{AC25} by studying $c$- and $g$-vectors under the assumption of sign-coherence. A fundamental difference from the classical cluster algebra setting is that sign-coherence does not hold in general for real-valued matrices. However, restricting to the sign-coherent case, they proved that finite types are classified by arbitrary finite Coxeter diagrams. This naturally generalizes the classical classification by Fomin and Zelevinsky \cite{FZ03a} over the integers.
	
	In this paper, we extend key structural theorems of finite type cluster algebras to the real-valued setting. Specifically, our primary objectives are to establish combinatorial characterizations of finite types and to prove the completeness of the corresponding $g$-fans within this generalized framework.
	
	There are two parts to our results. The first concerns weights alone. We call a weighted quiver \emph{$2$-finite} if every weight in its mutation class is less than $2$. Finite Coxeter quivers have this property, but for arbitrary real weights, a connected $2$-finite quiver need not be mutation-equivalent to a finite Coxeter quiver. We therefore consider the \emph{DB2-property}: every weight less than $2$ has the form
	\[
	2\cos\left(\frac{\pi}{m}\right)\qquad\text{for some }m\in\mathbb Z_{\geq3}.
	\]
	This property is motivated by the rank $2$ classification. In \cite{AC25}, they proved that every $c$-vector associated with the initial exchange quiver $Q$ of rank $2$ is sign-coherent if and only if $Q$ has the DB2-property. Hence, for any rank, the DB2-property is a necessary condition for sign-coherence. However, as demonstrated in Example~\ref{ex:DB2-not-SC}, this is not sufficient.
	\par
	The \emph{totally DB2-property} means that this condition holds throughout the mutation class. Under this condition, we obtain the following extension of the finite-type criterion in \cite{FZ03a}.
	
	\begin{maintheorem}[Theorem~\ref{thm:2fin-Cox}]\label{main:weights}
		Let $Q$ be a connected weighted quiver with the totally DB2-property. Then $Q$ is $2$-finite if and only if it is mutation-equivalent to a finite Coxeter quiver.
	\end{maintheorem}
	
	We say that a quiver is \emph{$g$-finite} if it has only finitely many $g$-clusters. Here, the totally SC-property means that, for every initial quiver in the mutation class, all $c$-vectors are sign-coherent and the $g$-vectors in each seed lie in a common closed orthant. Then, under the totally SC-property, we obtain the following characterization of $g$-finiteness.
	
	\begin{maintheorem}[Theorem~\ref{thm:gfin-Cox}]\label{main:finite}
		Let $Q$ be a weighted quiver with the totally SC-property. Then, the following conditions are equivalent.
		\begin{enumerate}[label=\textup{(\arabic*)}]
			\item $Q$ is $g$-finite; \label{item:gfinite}
			\item $Q$ is $2$-finite; \label{item:2finite}
			\item each connected component of $Q$ is mutation-equivalent to a finite Coxeter quiver. \label{item:coxeter}
		\end{enumerate}
	\end{maintheorem}
	
	The implication $\ref{item:coxeter}\Rightarrow\ref{item:gfinite}$ follows from \cite[Theorem~8.2]{AC25}, which shows that, under the totally SC-property, \ref{item:coxeter} holds if and only if every quiver in the mutation class of $Q$ is $g$-finite. The proof that this finiteness condition implies \ref{item:coxeter} uses the classification of mutation-finite quivers in \cite{FT23}. Here we deduce \ref{item:coxeter} from the $g$-finiteness of $Q$ alone by proving $\ref{item:gfinite}\Rightarrow\ref{item:2finite}$ directly in Proposition~\ref{prop:cfiniteto2finite} and then applying Theorem~\ref{main:weights} to obtain $\ref{item:2finite}\Rightarrow\ref{item:coxeter}$.
	
	The second part concerns the geometry of $g$-vectors and $g$-fans. For the construction of the $g$-fan and the formula for changing the initial seed, we use Assumption~\ref{ass:gfan}: the totally SC-property and the condition that, for every initial quiver in the mutation class, each $c$-vector on a coordinate axis is a signed standard basis vector.
	
	We say that a quiver is \emph{$g$-complete} if its $g$-fan covers the whole ambient space. For integer skew-symmetrizable exchange matrices, it was shown that finite type is equivalent to completeness of the $g$-fan \cite{N26,Y25}. We generalize this fact under Assumption~\ref{ass:gfan}.
	
	\begin{maintheorem}[Theorem~\ref{thm:main-classification}]
		\label{main:fans}
		Let $Q$ be a weighted quiver satisfying Assumption~\ref{ass:gfan}. Then the following conditions are equivalent:
		\begin{enumerate}[label=\textup{(\arabic*)}]
			\item $Q$ is $g$-finite;
			\item $Q$ is $g$-complete.
		\end{enumerate}
	\end{maintheorem}
	
	By Theorem~\ref{main:finite}, these conditions are also equivalent to each connected component of $Q$ being mutation-equivalent to a finite Coxeter quiver.
	
	The main step in the proof is the following restriction theorem. For a nonempty subset $I \subseteq Q_0$, write $Q_I$ for the full subquiver on $I$ and $\pi_I:\mathbb R^{Q_0}\to\mathbb R^I$ for the coordinate projection.
	
	\begin{maintheorem}[Proposition~\ref{prop:pi-in-cone} and Theorem~\ref{thm:I-gcomp-gdense}]
		\label{main:restriction}
		Suppose that $Q$ satisfies Assumption~\ref{ass:gfan}. For every $g$-cone $C$ of $Q$, there is a $g$-cone $\sigma$ of $Q_I$ such that $\pi_I(C)\subseteq\sigma$. Consequently, completeness and density of $g$-fans pass from $Q$ to $Q_I$.
	\end{maintheorem}
	
	For integer exchange matrices, the containment in Theorem~\ref{main:restriction} is proved in \cite[proof of Theorem~8]{CL20} as a corollary of \cite[Theorem~33]{Mu16}. Their original proof relied on the cluster scattering diagrams introduced by \cite{GHKK18}, and at present, this method is not applicable to our generalized setting. In contrast, we provide an alternative proof that relies on the initial quiver mutation formula introduced by \cite{NZ12}, which is cited as Lemma~\ref{lem:initial-seed-mut} in this paper. To prove Theorem~\ref{main:restriction}, we first show that the interior of a projected $g$-cone cannot meet a wall of the smaller fan. We then follow the projected cone along a sequence of mutations from the initial cone and show that at each step it stays in the same cone of the smaller fan or moves into an adjacent one. This argument does not assume that either fan is complete or locally finite. Together with mutation invariance of $g$-completeness and the rank-$2$ case, Theorem~\ref{main:restriction} applied to two vertices shows that a $g$-complete quiver cannot acquire an arrow of weight at least $2$ under mutation. Theorem~\ref{main:weights} then gives Theorem~\ref{main:fans}.
	
  The structure of this paper is as follows.
	Section~\ref{sec:weights} introduces $2$-finiteness and the DB2-property. Section~\ref{sec:sign} explains their relation with sign-coherence and $g$-finiteness. Section~\ref{sec:fans} proves the restriction theorem and the equivalence between $g$-finiteness and $g$-completeness. The weight classification is proved in Section~\ref{sec:pf-2fin-Cox}.

\section{Weighted quivers and \texorpdfstring{$2$}{2}-finiteness}\label{sec:weights}
\subsection{Weighted quivers and mutation}
	Throughout this paper, we fix a positive integer $n \in \mathbb{Z}_{\geq 1}$ called the \emph{rank}. A \emph{weighted quiver} is an $n$-vertex quiver without loops, $2$-cycles, or multiple arrows, in which each arrow is assigned a positive real weight. For brevity, we often refer to weighted quivers simply as quivers. Let $Q$ be a weighted quiver. We denote by $Q_0$ its vertex set and call $|Q_0|$ its rank. When the directions are irrelevant, we refer to arrows as edges. Unlabeled arrows and edges in diagrams have weight $1$. For $i,j\in Q_0$, put
	\[
	b_{ij}^Q=
	\begin{cases}
		w & \text{if there is an arrow $j\to i$ of weight $w$},\\
		-w & \text{if there is an arrow $i\to j$ of weight $w$},\\
		0 & \text{otherwise}.
	\end{cases}
	\]
	Then, the matrix $(b_{ij}^Q)_{i,j\in Q_0}$ is skew-symmetric.
	
	The \emph{mutation} of $Q$ at $k\in Q_0$ is a weighted quiver $\mu_k(Q)$ defined by
	\begin{equation}\label{eq:quiver-mutation}
		b_{ij}^{\mu_k(Q)}=
		\begin{cases}
			-b_{ij}^Q & \text{if $i=k$ or $j=k$},\\
			b_{ij}^Q+[b_{ik}^Q]_+[b_{kj}^Q]_+
			-[-b_{ik}^Q]_+[-b_{kj}^Q]_+ & \text{otherwise},
		\end{cases}
	\end{equation}
	where $[a]_+:=\max\{a,0\}$ for $a\in\mathbb R$. Mutation is an involution, that is $\mu_k(\mu_k(Q))=Q$. Two weighted quivers are \emph{mutation-equivalent} if they are related by a finite sequence of mutations, up to relabeling of vertices. The \emph{mutation class} of $Q$ is the set of weighted quivers mutation-equivalent to $Q$. Throughout the paper, full subquivers are understood to have nonempty vertex sets. For a nonempty subset $I\subseteq Q_0$, we denote by $Q_I$ the full subquiver of $Q$ on the vertex set $I$, with the weights inherited from $Q$. When $I = Q_0 \setminus \{k\}$, we write $Q \setminus \{k\}$ for $Q_I$. The mutation formula gives
	\begin{equation}\label{eq:restriction-mutation}
		(\mu_k(Q))_I=\mu_k(Q_I)\qquad(k\in I).
	\end{equation}

\subsection{\texorpdfstring{$2$}{2}-finiteness and the discrete-below-$2$ property}\label{sec:2fin-DB2}
	We introduce two conditions on the weights in a mutation class and show that together they characterize finite Coxeter type.
	
	\begin{definition}\label{def:2finite}
		A weighted quiver is called \emph{$2$-finite} if every weight in every quiver in its mutation class is less than $2$, and \emph{$2$-infinite} otherwise.
	\end{definition}
	
	We next introduce a condition on the weights below $2$. For $m\in\mathbb Z_{\geq2}$, set
	\[
	\bra{m}:=2\cos\left(\frac{\pi}{m}\right),
	\qquad
	\mathcal W_{<2}:=\{\bra{m}\mid m\in\mathbb Z_{\geq2}\}.
	\]
	For example, $\bra{2}=0$, $\bra{3}=1$, $\bra{4}=\sqrt2$, and $\bra{6}=\sqrt3$. In particular, $\bra{m}\geq1$ for $m\in\mathbb Z_{\geq3}$.
	
	\begin{definition}\label{def:DB2}
		A weighted quiver has the \emph{discrete-below-$2$ property} (\emph{DB2-property}) if each of its weights belongs to $\mathcal W_{<2}\sqcup\mathbb R_{\geq2}$. It has the \emph{totally DB2-property} if every quiver in its mutation class has the DB2-property.
	\end{definition}
	
	We denote by $\D$ the set of weighted quivers with the totally DB2-property. If $Q\in\D$ is $2$-finite, then every weight in its mutation class is $\bra{m}$ for some integer $m\geq3$, and is therefore at least $1$. We record three basic facts about $2$-finiteness and the totally DB2-property.
	\begin{enumerate}
		\item Both properties are invariant under mutation.
		\item Both properties pass to full subquivers, since every sequence of mutations on $Q_I$ extends to the same sequence on $Q$, by~\eqref{eq:restriction-mutation}.
		\item Each property holds for $Q$ if and only if it holds for every connected component of $Q$, since mutations on distinct components are independent.
	\end{enumerate}
	
	We now introduce a class of quivers that are $2$-finite and have the totally DB2-property.
	
\begin{figure}[htbp]
		\centering
		  %%%%%%%%%%%% type A_n
		\begin{tikzpicture}[baseline = {(0,0)}]
			\foreach \x in {1,2,3,4,5,6,7}
			{
				\node[myQV] (\x) at (\x,0) {};
			}
			\foreach \x [evaluate = \x as \y using int(\x+1)] in {1,2,3,5,6}
			{
				\draw[-] (\x) -> (\y);
			}
			\draw[dashed] (4) -- (5);
			\node[anchor = east] at (1.west) [left] {$A_n$};
		\end{tikzpicture}
		\\
		%%%%%%%%%%%% type B_n
		\begin{tikzpicture}[baseline = {(0,0)}]
			\foreach \x in {1,2,3,4,5,6,7}
			{
				\node[myQV] (\x) at (\x,0) {};
			}
			\foreach \x [evaluate = \x as \y using int(\x+1)] in {1,2,3,5}
			{
				\draw[-] (\x) -> (\y);
			}
			\draw[dashed] (4) -- (5);
			\draw[-] (6) -> (7) node [pos = 0.5, above] {$\bra{4}$};
			\node[anchor = east] at (1.west) [left] {$B_n$};
		\end{tikzpicture}
		\\
		%%%%%%%%%%%% type D_n
		\begin{tikzpicture}[baseline = {(0,0)}]
			\foreach \x in {1,2,3,4,5,6}
			{
				\node[myQV] (\x) at (\x,0) {};
			}
			\node[myQV] (7) at ({6 + cos(30)}, {sin(30)}) {};
			\node[myQV] (8) at ({6 + cos(-30)}, {sin(-30)}) {};
			\foreach \x [evaluate = \x as \y using int(\x+1)] in {1,2,3,5,6}
			{
				\draw[-] (\x) -> (\y);
			}
			\draw[dashed] (4) -- (5);
			\draw[-] (6) -> (8);
			\node[anchor = east] at (1.west) [left] {$D_n$};
		\end{tikzpicture}
		\\
		%%%%%%%%%%%%% type E_6
		\begin{tikzpicture}[baseline = {(0,0)}]
			\node[myQV] (1) at (1,0) {};
			\foreach \x [evaluate = \x as \y using int(\x-1)] in {2,3,4,5}
			{
				\node[myQV] (\x) at (\x,0) {};
				\draw[-] (\y) -> (\x);
			}
			\node[myQV] (3') at (3,1) {};
			\draw[-] (3) -> (3');
			\node[anchor = east] at (1.west) [left] {$E_6$};
		\end{tikzpicture}
		\\
		%%%%%%%%%%%%% type E_7
		\begin{tikzpicture}[baseline = {(0,0)}]
			\node[myQV] (1) at (1,0) {};
			\foreach \x [evaluate = \x as \y using int(\x-1)] in {2,3,4,5,6}
			{
				\node[myQV] (\x) at (\x,0) {};
				\draw[-] (\y) -> (\x);
			}
			\node[myQV] (3') at (3,1) {};
			\draw[-] (3) -> (3');
			\node[anchor = east] at (1.west) [left] {$E_7$};
		\end{tikzpicture}
		\\
		%%%%%%%%%%%%% type E_8
		\begin{tikzpicture}[baseline = {(0,0)}]
			\node[myQV] (1) at (1,0) {};
			\foreach \x [evaluate = \x as \y using int(\x-1)] in {2,3,4,5,6,7}
			{
				\node[myQV] (\x) at (\x,0) {};
				\draw[-] (\y) -> (\x);
			}
			\node[myQV] (3') at (3,1) {};
			\draw[-] (3) -> (3');
			\node[anchor = east] at (1.west) [left] {$E_8$};
		\end{tikzpicture}
		\\
		%%%%%%%%%%%%% type F_4
		\begin{tikzpicture}
			\foreach \x in {1,2,3,4}
			{
				\node[myQV] (\x) at (\x,0) {};
			}
			\draw[-] (1) -> (2);
			\draw[-] (2) -> (3) node [pos = 0.5, above] {$\bra{4}$};
			\draw[-] (3) -> (4);
			\node[anchor = east] at (1.west) [left] {$F_4$};
		\end{tikzpicture}
		\\
		%%%%%%%%%%%%% type H_3
		\begin{tikzpicture}
			\foreach \x in {1,2,3}
			{
				\node[myQV] (\x) at (\x,0) {};
			}
			\draw[-] (1) -> (2);
			\draw[-] (2) -> (3) node [pos = 0.5, above] {$\bra{5}$};
			\node[anchor = east] at (1.west) [left] {$H_3$};
		\end{tikzpicture}
		\\
		%%%%%%%%%%%%% type H_4
		\begin{tikzpicture}
			\foreach \x in {1,2,3,4}
			{
				\node[myQV] (\x) at (\x,0) {};
			}
			\draw[-] (1) -> (2);
			\draw[-] (2) -> (3);
			\draw[-] (3) -> (4) node [pos = 0.5, above] {$\bra{5}$};
			\node[anchor = east] at (1.west) [left] {$H_4$};
		\end{tikzpicture}
		\\
		%%%%%%%%%%%%% type I_2(m)
		\begin{tikzpicture}
			\node[myQV] (1) at (1,0) {};
			\node[myQV] (2) at (2,0) {};
			\draw[-] (1) -> (2) node [pos = 0.5, above] {$\bra{m}$};
			\node[anchor = east] at (1.west) [left] {$I_2(m)$};
		\end{tikzpicture}
		\caption{Underlying weighted graphs of finite Coxeter quivers}
		\label{fig:finCox}
	\end{figure}
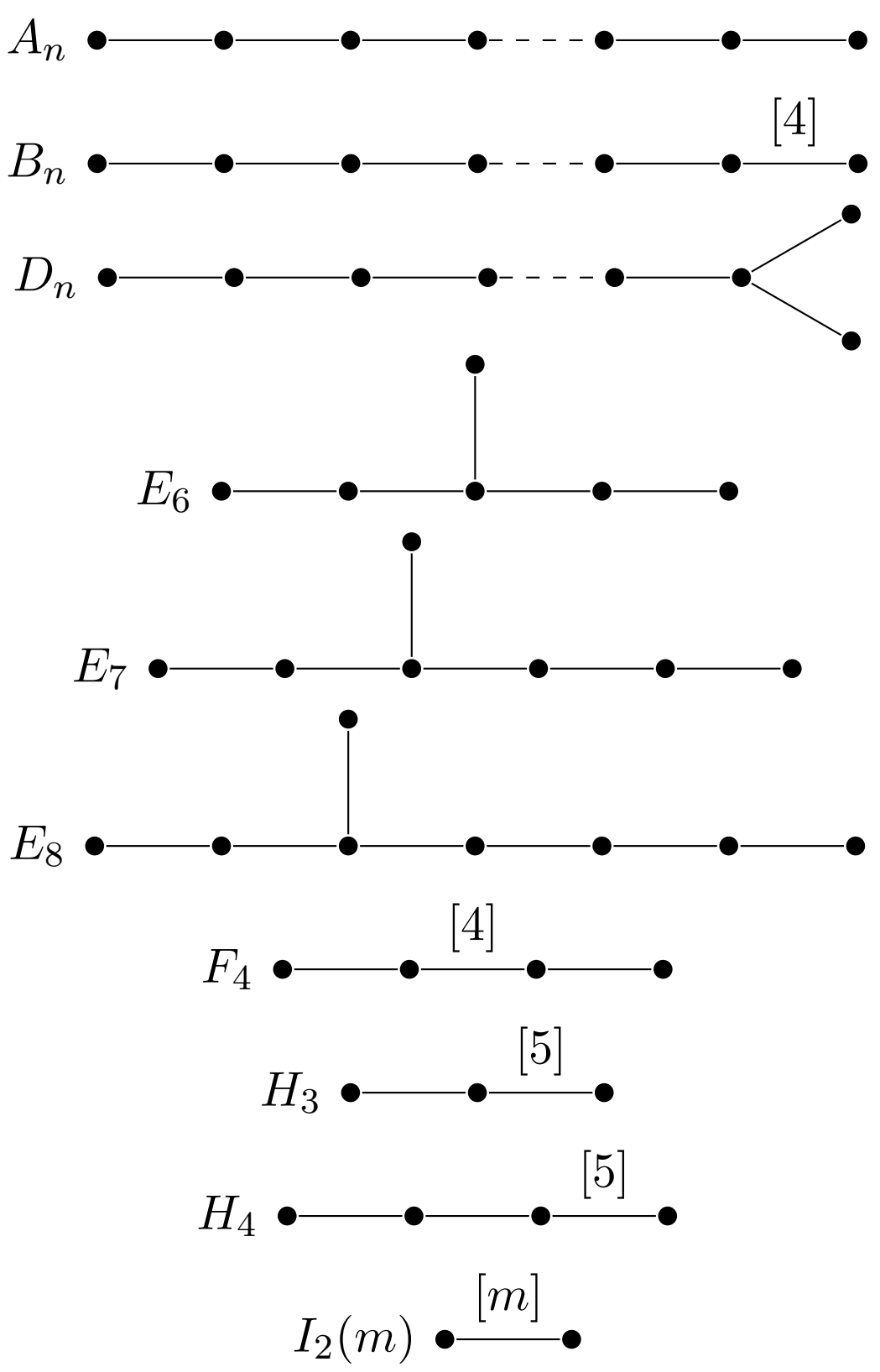

	A \emph{finite Coxeter quiver} is a weighted quiver obtained by orienting the Coxeter diagram of an irreducible finite Coxeter group and replacing each edge label $m$ by the weight $\bra{m}$; see Figure~\ref{fig:finCox}. The subscript in each type denotes the number of vertices. A weighted quiver is said to be \emph{of finite Coxeter type} if each of its connected components is mutation-equivalent to a finite Coxeter quiver. Note that the one-vertex quiver is included as type $A_1$, and $I_2(3)$, $I_2(4)$, and $I_2(6)$ are of types $A_2$, $B_2$, and $G_2$, respectively. Any two orientations of a weighted tree are related by source and sink mutations, so the orientation does not affect the mutation class; see Lemma~\ref{lem: chain reorientations}.
	
	By \cite[Theorem~8.2 and Proposition~8.5(b)]{AC25} (see also Theorem~\ref{thm:Cox-SC-gfin}), every weighted quiver $Q$ satisfies
	\begin{equation}\label{eq:Cox-2fin-DB2}
		Q\text{ is of finite Coxeter type}
		\quad\Longrightarrow\quad
		Q\text{ is $2$-finite and }Q\in\D.
	\end{equation}
	The disconnected case follows by considering connected components. The following example shows that $2$-finiteness alone does not imply finite Coxeter type.
	
	\begin{example}\label{ex:2fin-non-Cox}
		Let $Q$ be the directed triangle with all three weights equal to $1/2$. Mutating at any vertex gives an acyclic triangle with weights $1/2,1/2,1/4$. Since the Markov constant of $Q$ is $3(1/2)^2-(1/2)^3=5/8$, it follows from \cite[Proposition~1.2]{CK26} that every weight in its mutation class is at most $\sqrt{5/8}<2$. Thus $Q$ is $2$-finite. Since $1/2\notin\mathcal W_{<2}$, it does not have the DB2-property, and hence it is not of finite Coxeter type by \eqref{eq:Cox-2fin-DB2}.
	\end{example}
	
	Our first main result shows that $2$-finiteness together with the totally DB2-property characterizes finite Coxeter type.
	
	\begin{theorem}\label{thm:2fin-Cox}
		A weighted quiver is $2$-finite and has the totally DB2-property if and only if it is of finite Coxeter type.
	\end{theorem}
	
	The ``if'' implication is~\eqref{eq:Cox-2fin-DB2}. The converse will be proved in Section~\ref{sec:pf-2fin-Cox}.

\section{Sign-coherence and finiteness}\label{sec:sign}
	Throughout this section, let $Q$ be a weighted quiver with $Q_0=\{1,\ldots,n\}$.
\subsection{\texorpdfstring{$g$}{g}-vectors and \texorpdfstring{$c$}{c}-vectors}
	The \emph{framed quiver} of $Q$ is a weighted quiver $\hat Q$ obtained from $Q$ by adding vertices $n+1,\ldots,2n$ and arrows $i\to n+i$ of weight $1$ for $1\leq i\leq n$. Mutations are allowed only at the original vertices $1,\ldots,n$. We will refer to the full subquiver on these original vertices as the \emph{mutable part}.

  Let $\mathbf e_1,\ldots,\mathbf e_n$ be the standard basis of $\mathbb R^n$, and
	\[
	\mathbf b_k^Q:=(b_{1k}^Q,\ldots,b_{nk}^Q)^{\mathsf T}
	\qquad(1\leq k\leq n).
	\]
	
	\begin{definition}\label{def:gseed}
		A \emph{$g$-vector seed for $Q$} is a pair $(R,(\mathbf g_1,\ldots,\mathbf g_n))$ obtained from the initial seed
		\[
		(\hat Q,(\mathbf e_1,\ldots,\mathbf e_n))
		\]
		by a finite sequence of mutations. Here, mutation of $(R,(\mathbf g_1,\ldots,\mathbf g_n))$ at $k$ replaces $R$ by $\mu_k(R)$ and replaces $(\mathbf g_1,\ldots,\mathbf g_n)$ by $(\mathbf g_1,\ldots, \mathbf{g}'_k, \ldots ,\mathbf g_n)$, where
		\begin{equation}\label{eq:g-final-mutation}
			\mathbf g'_k
      =
      -\mathbf g_k + \sum_{i=1}^n[b_{ik}^R]_+\mathbf g_i - \sum_{j=1}^n[b_{n+j,k}^R]_+\mathbf b_j^Q.
		\end{equation}
		The unordered set $\{\mathbf g_1,\ldots,\mathbf g_n\}$ is a \emph{$g$-cluster} for $Q$. Its elements are \emph{$g$-vectors for $Q$}. The \emph{$c$-vectors} of the $g$-vector seed are the vectors
		\[
		\mathbf c_j^R:=(b_{n+1,j}^R,\ldots,b_{2n,j}^R)^{\mathsf T}
		\qquad(1\leq j\leq n).
		\]
	\end{definition}
	
	The matrices $G=(\mathbf g_1\ \cdots\ \mathbf g_n)$ and $C=(\mathbf c_1^R\ \cdots\ \mathbf c_n^R)$ are called the \emph{$G$-matrix} and the \emph{$C$-matrix} of the seed, respectively. Both initial matrices are the identity matrix $I_n$. We denote by $\gseed_n(Q)$ the set of $g$-clusters for $Q$ and, for $0\leq k\leq n$, put
	\[
	\gseed_k(Q):=
	\left\{g\subseteq\mathbb R^n\relmiddle|
	\#g=k,\quad g\subseteq g^\ast\text{ for some }g^\ast\in\gseed_n(Q)\right\}.
	\]
	In particular, $\gseed_0(Q)=\{\varnothing\}$.
	
	\begin{definition}\label{def:sign-coherence}
		A vector is \emph{sign-coherent} if all its coordinates are nonnegative or all are nonpositive. A $g$-cluster $\{\mathbf g_1,\ldots,\mathbf g_n\}$ is \emph{sign-coherent} if, for each $i$, the numbers $(\mathbf g_1)_i,\ldots,(\mathbf g_n)_i$ are all nonnegative or all nonpositive. Equivalently, the vectors in the seed lie in a common closed orthant.
	\end{definition}
	
	\begin{definition}\label{def:SC}
		A weighted quiver has the \emph{sign-coherent property} (\emph{SC-property}, for short) if all its $c$-vectors and all its $g$-clusters are sign-coherent. It has the \emph{totally SC-property} if every quiver in its mutation class has the SC-property.
	\end{definition}
	
	It is conjectured that the SC-property implies the totally SC-property
	\cite[Conjecture~5.1]{AC25}; see also \cite[Section~8]{Re14}.
	
	\begin{definition}\label{def:gfinite}
		A weighted quiver $Q$ is called \emph{$g$-finite} if $\gseed_n(Q)$ is finite.
	\end{definition}
	
	Suppose that $Q$ has the SC-property. Consider a $g$-vector seed $(R,(\mathbf{g}_1,\dots,\mathbf{g}_n))$ for $Q$ and the corresponding $c$-vectors $(\mathbf{c}_1,\dots,\mathbf{c}_n)$. For a vertex $k=1,\dots,n$, let $\varepsilon_k \in \{1,-1\}$ be the sign of $\mathbf{c}_{k}$. Then the mutation at $k$ of the $g$-vectors $(\mathbf{g}_1,\dots,\mathbf{g}_n) \mapsto (\mathbf{g}_1,\dots,\mathbf{g}'_k,\dots,\mathbf{g}_n)$ and the $c$-vectors $(\mathbf{c}_1,\dots,\mathbf{c}_n) \mapsto (\mathbf{c}'_1,\dots,\mathbf{c}'_n)$ is given by the following basis transformations:
  \begin{align}
    \mathbf{g}'_{k} &= -\mathbf{g}_k + \sum_{i=1}^n[-\varepsilon_k b_{ik}^R]_+\mathbf g_i, \label{eq:g-signed-mutation}
    \\
    \mathbf{c}'_{\ell} &= \begin{cases}
      -\mathbf{c}_k & \textup{if $\ell = k$},
      \\
      \mathbf{c}_{\ell} + [\varepsilon_k b_{k\ell}^R]_+\mathbf{c}_{k} & \textup{if $\ell \neq k$};
    \end{cases}
  \end{align}
  see \cite{AC25}. Thus every $g$-cluster is a basis of $\mathbb R^n$.
  Moreover the second duality in \cite[Proposition~4.8]{AC25} gives \begin{equation}\label{eq:GC-duality}
		G^{\mathsf T}C=I_n,
	\end{equation}
	and hence $C=(G^{-1})^{\mathsf T}$ and $G = (C^{-1})^{\mathsf T}$ hold. Since an unordered $g$-cluster has at most $n!$ orderings, the following conditions are equivalent:
	\begin{equation}\label{eq:equivalenceofcgfinite}
    \begin{gathered}
    Q\text{ is $g$-finite}
    \quad\Longleftrightarrow\quad
    \text{there are finitely many $G$-matrices}\\
    \quad\Longleftrightarrow\quad
    \text{there are finitely many $C$-matrices}
    \quad\Longleftrightarrow\quad
    \text{there are finitely many $c$-vectors}.
    \end{gathered}
  \end{equation}
	We use these equivalences without further mention and state results about the last condition from \cite{AC25} in terms of $g$-finiteness whenever the SC-property holds. We recall the following result.
	
	\begin{theorem}[{\cite[Theorem~8.2]{AC25}}]\label{thm:Cox-SC-gfin}
		Every weighted quiver of finite Coxeter type has the totally SC-property and is $g$-finite.
	\end{theorem}

\subsection{Sign-coherence and the DB2-property}
  Here, we collect some basic observations concerning sign-coherence and the DB2-property.

  In the rank $2$ case, sign-coherence is completely classified as follows:
	\begin{theorem}[{\cite[Theorem~7.1 and Example~7.5]{AC25}}]\label{thm:SC-D-rank2}
		Let $Q$ be a weighted quiver of rank $2$. Then $Q$ has the SC-property if and only if it has the DB2-property. Under these equivalent conditions, $Q$ is $g$-finite if and only if it is $2$-finite.
	\end{theorem}
	
  Note that the above classification provides a necessary condition for sign-coherence in arbitrary rank via the following lemma.

	\begin{lemma}\label{lem:SC-full}
		For any nonempty subset $I\subseteq Q_0$, we have the following claims.
		\begin{enumerate}[label=\textup{(\arabic*)}]
			\item If $Q$ has the SC-property, then $Q_I$ has the SC-property.
			\item If $Q$ has the totally SC-property, then $Q_I$ has the totally SC-property.
			\item If $Q$ is $g$-finite, then $Q_I$ is $g$-finite.
		\end{enumerate}
	\end{lemma}
	
	\begin{proof}
		Start with the framed quivers of $Q$ and $Q_I$, and apply the same sequence of mutations at vertices in $I$. Denote by $R$ and $S$ the resulting quivers, respectively. Let $\iota_I:\mathbb R^I\to\mathbb R^{Q_0}$ be the coordinate inclusion and $\pi_I:\mathbb R^{Q_0}\to\mathbb R^I$ the coordinate projection. Since only vertices in $I$ are mutated, we have $\mathbf{g}_k=\mathbf{e}_k$ for $k\notin I$. Moreover, the coordinates outside $I$ of $\mathbf{c}_j^R$ remain zero for $j\in I$. Using these facts and the mutation formulas \eqref{eq:quiver-mutation} and \eqref{eq:g-final-mutation}, the identities
		\[
		(b_{ij}^R)_{i,j\in I}=(b_{ij}^S)_{i,j\in I},\qquad
		\mathbf c_j^R=\iota_I(\mathbf c_j^S),\qquad
		\pi_I(\mathbf g_j^R)=\mathbf g_j^S\quad(j\in I)
		\]
		follow by induction on the length of the sequence of mutations.
		
		By the second and third identities, the columns of the smaller $C$-matrix and the rows of its $G$-matrix inherit sign-coherence, which proves (1). The third identity also proves (3).
		
		For (2), suppose that $Q$ has the totally SC-property. Let $\mu$ be any sequence of mutations at vertices in $I$. By \eqref{eq:restriction-mutation}, we have $\mu(Q_I)=(\mu(Q))_I$. Since $\mu(Q)$ has the SC-property by assumption, (1) implies that $\mu(Q_I)$ also has the SC-property. Hence $Q_I$ has the totally SC-property.
	\end{proof}
	
	\begin{corollary}\label{cor:SC-D}
		The SC-property implies the DB2-property. The totally SC-property implies the totally DB2-property.
	\end{corollary}
	
	\begin{proof}
		The first assertion follows by applying Lemma~\ref{lem:SC-full}(1) and Theorem~\ref{thm:SC-D-rank2} to the full subquiver on the endpoints of each arrow. Applying the first assertion to every quiver in the mutation class gives the second.
	\end{proof}
	
	The converse of the first assertion fails in general.
	
	\begin{example}\label{ex:DB2-not-SC}
		Consider
		\[
		Q=\begin{tikzpicture}[baseline={(0,0)},scale=0.75]
			\node[myCN] (1) at (90:1) {$1$};
			\node[myCN] (2) at (210:1) {$2$};
			\node[myCN] (3) at (330:1) {$3$};
			\draw[->] (1) -> node[pos=0.5,auto=right] {$\bra4$} (2);
			\draw[->] (3) -> (1);
			\draw[->] (2) -> (3);
		\end{tikzpicture}.
		\]
		Then $Q$ has the DB2-property. After mutations at $1$, $2$, and $3$, in that order, the resulting $C$- and $G$-matrices are
		\[
		C=
		\begin{bmatrix}
			\sqrt2-2&\sqrt2-1&-1\\
			3-2\sqrt2&2-2\sqrt2&1-\sqrt2\\
			\sqrt2-1&\sqrt2-1&-1
		\end{bmatrix},\qquad
		G=
		\begin{bmatrix}
			-1&0&1-\sqrt2\\
			0&-1&1-\sqrt2\\
			1&\sqrt2-1&1-\sqrt2
		\end{bmatrix},
		\]
		where $\bra4=\sqrt2$. Since $1<\sqrt2<2$, each of the first two columns of $C$ has both positive and negative entries, and the third row of $G$ also has entries of both signs. Thus $Q$ does not have the SC-property. On the other hand, mutation at $1$ already produces a weight $\sqrt2-1$, which does not belong to $\mathcal W_{<2}\sqcup\mathbb R_{\geq2}$. Hence this example does not have the totally DB2-property.
	\end{example}
	
	\begin{corollary}\label{cor:2fin-TSC-TD}
		Let $Q$ be $2$-finite. Then $Q$ has the totally SC-property if and only if it has the totally DB2-property.
	\end{corollary}
	
	\begin{proof}
		The ``only if'' part follows from Corollary~\ref{cor:SC-D}. Conversely, if $Q$ has the totally DB2-property, then it is of finite Coxeter type by Theorem~\ref{thm:2fin-Cox} and hence has the totally SC-property by Theorem~\ref{thm:Cox-SC-gfin}.
	\end{proof}
	
	Theorem~\ref{thm:SC-D-rank2} and Corollary~\ref{cor:2fin-TSC-TD} suggest the following converse to the second assertion in Corollary~\ref{cor:SC-D}.
	
	\begin{conjecture}\label{conj:TD-TSC}
		Every weighted quiver with the totally DB2-property has the totally SC-property.
	\end{conjecture}
  In \cite[Conjecture~5.1]{AC25}, they conjectured that the SC-property implies the totally SC-property. By combining these two conjectures, we expect that the SC-property, the totally SC-property, and the totally DB2-property are all equivalent.

\subsection{Finiteness and the Coxeter classification}	
	We next characterize $g$-finiteness under the totally SC-property. The main step is to show the following claim directly.
  \begin{proposition}\label{prop:cfiniteto2finite}
    Let $Q$ be a weighted quiver with the SC-property. If $Q$ is $g$-finite, then $Q$ is $2$-finite.
  \end{proposition}
  Note that in the following proof, we do not choose a particular initial quiver $Q$ from its mutation class. In particular, when $Q$ is $2$-infinite and has the SC-property, we can conclude that $Q$ is not $g$-finite even if all weights of $Q$ itself are less than $2$.
  \begin{proof}
    By \eqref{eq:equivalenceofcgfinite}, there are only finitely many $c$-vectors.
    Suppose, for a contradiction, that a quiver $R$ in the mutation class of $Q$ has an arrow of weight $s \geq 2$.  Choose two vertices $i,j$ in $R$ such that $b_{ij}^R = s$. Denote the corresponding pair of $c$-vectors by $(\mathbf c_i,\mathbf c_j)$.
		
		Starting from this quiver, mutate alternately at $i$ and $j$. Set $\mathbf v_{-1}=-\mathbf c_j$ and $\mathbf v_0=\mathbf c_i$. By \eqref{eq:quiver-mutation}, the successive pairs of $c$-vectors are
    \[
		(\mathbf v_0,-\mathbf v_{-1}) \overset{\mu_i}{\mapsto}
		(-\mathbf v_0,\mathbf v_1) \overset{\mu_j}{\mapsto}
		(\mathbf v_2,-\mathbf v_1) \overset{\mu_i}{\mapsto}
		(-\mathbf v_2,\mathbf v_3) \overset{\mu_j}{\mapsto} \cdots
		\]
    where these vectors $\mathbf{v}_r$ are recursively obtained by
		\begin{equation}\label{eq:recursionofv}
		  \mathbf v_{r+1}=s[\mathbf v_r]_+-\mathbf v_{r-1}
		  \qquad(r\geq0).
    \end{equation}
		Here, $[\mathbf v]_+ = ([v_1]_{+},\dots,[v_n]_{+})^{\mathsf T}$ for $\mathbf{v}=(v_1,\dots,v_n)^{\mathsf T} \in \mathbb{R}^n$. 
    By the assumption, the vectors $\mathbf{v}_r$ take only finitely many values. Since each mutation is an involution, this means that the sequence $(\mathbf{v}_r)$ must be periodic.

		On the other hand, by $\mathbf{v}_r = [\mathbf{v}_r]_{+} - [-\mathbf{v}_r]_{+}$ and \eqref{eq:recursionofv}, we have the coordinatewise inequality
		\[
		(\mathbf{v}_{r+1}-\mathbf{v}_r)-(\mathbf{v}_r-\mathbf{v}_{r-1})
		=(s-2)[\mathbf v_r]_++2[-\mathbf v_r]_+
		\geq 0.
		\]
		Thus, in each coordinate, the successive differences $\mathbf{v}_{r+1}-\mathbf{v}_{r}$ are nondecreasing. Since these differences are periodic, $\mathbf{v}_{r+1} - \mathbf{v}_r$ must be constant. Periodicity of $(\mathbf v_r)$ then forces this constant to be zero, and hence $\mathbf{v}_r$ is independent of $r$. In particular, $\mathbf c_i=\mathbf v_0=\mathbf v_{-1}=-\mathbf c_j$, contradicting the invertibility of the $C$-matrix. Thus $Q$ is $2$-finite.
  \end{proof}
	In conclusion, we obtain the following characterization of $g$-finiteness.
	\begin{theorem}\label{thm:gfin-Cox}
		Let $Q$ be a weighted quiver with the totally SC-property. Then the following conditions are equivalent:
		\begin{enumerate}[label=\textup{(\arabic*)}]
			\item $Q$ is $g$-finite;
			\item $Q$ is $2$-finite;
			\item $Q$ is of finite Coxeter type.
		\end{enumerate}
	\end{theorem}
	
	\begin{proof}
		By Corollary~\ref{cor:SC-D}, $Q$ has the totally DB2-property. Hence (2)$\Leftrightarrow$(3) follows from Theorem~\ref{thm:2fin-Cox}, and (3)$\Rightarrow$(1) follows from Theorem~\ref{thm:Cox-SC-gfin}. The remaining implication (1)$\Rightarrow$(2) follows from Proposition~\ref{prop:cfiniteto2finite}.
	\end{proof}
	
	Since $2$-finiteness is mutation-invariant, Theorem~\ref{thm:gfin-Cox} also shows that $g$-finiteness is invariant under mutation among weighted quivers with the totally SC-property.

\section{Completeness and density of \texorpdfstring{$g$}{g}-fans}\label{sec:fans}
	Throughout this section, let $Q$ be a weighted quiver with $Q_0=\{1,\ldots,n\}$. Under the following assumption, we introduce the $g$-fan of $Q$ and prove that it is complete if and only if $Q$ is $g$-finite.
	
	\begin{assumption}\label{ass:gfan}
		The weighted quiver $Q$ satisfies the following conditions:
		\begin{enumerate}[label=\textup{(\arabic*)}]
			\item $Q$ has the totally SC-property;
			\item for every quiver $Q'$ in the mutation class of $Q$, every $c$-vector for $Q'$ lying on a coordinate axis is equal to $\pm\mathbf e_i$ for some $i\in Q'_0$.
		\end{enumerate}
	\end{assumption}
	
	Unless stated otherwise, Assumption~\ref{ass:gfan} is in force throughout this section. Condition~(2) is the coordinate-axis condition in \cite[Conjecture~5.3]{AC25}, imposed throughout the mutation class. For skew-symmetric exchange matrices, Assumption~\ref{ass:gfan} is the hypothesis used in \cite[Conjecture~5.10 and Section~6]{AC25}. It holds for quivers of quasi-integer type, for rank-$2$ quivers with the SC-property, and for weighted quivers of finite Coxeter type \cite[Proposition~5.11, Theorems~7.7 and~8.9]{AC25}.

\subsection{The $g$-fan}
	\begin{lemma}\label{lem:assumption-restriction}
		Assumption~\ref{ass:gfan} is preserved by mutation and passes to full subquivers.
	\end{lemma}
	
	\begin{proof}
		Mutation invariance is part of the assumption. Let $I\subseteq Q_0$. The totally SC-property passes to $Q_I$ by Lemma~\ref{lem:SC-full}(2). To check condition (2), let $S$ be any quiver in the mutation class of $Q_I$. Lift the sequence of mutations to $Q$, obtaining a quiver $R$ with $R_I=S$. Use $R$ and $S$ as the initial quivers for the $c$-vectors. The proof of Lemma~\ref{lem:SC-full} shows that every $c$-vector for $S$ lifts, by adding zero coordinates outside $I$, to a $c$-vector for $R$ along the same sequence of mutations. If a $c$-vector for $S$ lies on a coordinate axis, then its lift lies on the same axis and is therefore equal to $\pm\mathbf e_i$ for some $i\in I$ by condition~(2) for $R$. This proves condition~(2) for $S$.
	\end{proof}
	
	For a finite set $g=\{\mathbf g_1,\ldots,\mathbf g_k\}$ of vectors, a \emph{(polyhedral) cone} is a closed set of the form
	\[
	C(g)=C(\mathbf g_1,\ldots,\mathbf g_k)
	:=\left\{\sum_{i=1}^k a_i\mathbf g_i\relmiddle|a_i\geq0\right\}.
	\]
	We use $C(\varnothing)=\{\mathbf0\}$. A \emph{face} of a cone $C\subseteq\mathbb R^n$ is a set $C\cap\ker\ell$, where $\ell:\mathbb R^n\to\mathbb R$ is a linear functional that is nonnegative on $C$. A cone is called \emph{simplicial} if it is generated by a linearly independent set. If $g$ is linearly independent, then the faces of $C(g)$ are precisely the cones $C(g')$ with $g'\subseteq g$.
  
  	A cone $C(g)$ with $g\in\gseed_n(Q)$ is called a \emph{$g$-cone}. For $n\geq1$, a cone $C(g)$ with $g\in\gseed_{n-1}(Q)$ is called a \emph{wall}. Set
	\[
	\mathcal F(Q):=\{C(g)\mid g\in\gseed_k(Q),\ 0\leq k\leq n\}.
	\]
	A \emph{fan} is a collection of polyhedral cones closed under taking faces, such that the intersection of any two cones is a face of each. A fan is called \emph{simplicial} if all its cones are simplicial. We do not require fans to be finite or locally finite.
	
	\begin{proposition}[{\cite[Proposition~6.6]{AC25}}]\label{prop:fan}
		Let $Q$ satisfy Assumption~\ref{ass:gfan}. Then $\mathcal F(Q)$ is a simplicial fan.
	\end{proposition}
	
	\begin{proof}
		The fan property is proved in \cite[Proposition~6.6]{AC25}. Every $g$-cluster is a basis by \eqref{eq:GC-duality}, so its cones are simplicial.
	\end{proof}
	
	We call $\mathcal F(Q)$ the \emph{$g$-fan} of $Q$. Thus the term $g$-cone refers to a full-dimensional cone of this fan. The support of the fan is
	\[
	|\mathcal F(Q)|=\bigcup_{g\in\gseed_n(Q)}C(g).
	\]
	Two distinct $g$-cones are \emph{adjacent} if they share a wall. We use $\int$ for interior in the ambient vector space. Relative interiors are taken in the linear span of the cone.
	
	\begin{definition}\label{def:complete-dense}
		Let $Q$ satisfy Assumption~\ref{ass:gfan}. It is called \emph{$g$-complete} if $|\mathcal F(Q)|=\mathbb R^{n}$, and \emph{$g$-dense} if $|\mathcal F(Q)|$ is dense in $\mathbb R^{n}$.
	\end{definition}
	
	To prove that $g$-finiteness implies $g$-completeness, we first record the following basic properties of $g$-cones.
	
	\begin{lemma}\label{lem:two-cones}
		Let $Q$ satisfy Assumption~\ref{ass:gfan}. Then, every wall lies in exactly two $g$-cones. These $g$-cones lie on opposite sides of the hyperplane spanned by the wall, and their union contains a neighborhood of every point in the relative interior of the wall. Moreover, any two $g$-cones can be joined by a finite sequence of adjacent $g$-cones.
	\end{lemma}
	
	\begin{proof}
		For a given wall $W$, choose a $g$-vector seed $(R,(\mathbf g_1,\ldots,\mathbf g_n))$ such that the wall is generated by all $g$-vectors except $\mathbf g_k$. Consider the mutation of $g$-vectors at $k$. Let $\mathbf{g}'_k$ be the generated $g$-vector given by \eqref{eq:g-signed-mutation}. Since $b_{kk}^R = 0$, we obtain $\mathbf g'_k = -\mathbf g_k$ in $\mathbb{R}^{n}/\langle W \rangle$, where $\langle W \rangle \subset \mathbb{R}^{n}$ is the linear subspace spanned by $W$. Thus mutation at $k$ gives a second $g$-cone sharing the same wall $W$ on the opposite side. Since both cones are full-dimensional and simplicial, each cone fills one of the two half-neighborhoods near a point in the relative interior of their common wall. Their union therefore contains a neighborhood of that point.
		
		There cannot be a third full-dimensional cone containing the same wall: on either side, it would have interior intersecting one of these two cones, contrary to Proposition~\ref{prop:fan}. Finally, every $g$-vector seed is obtained from the initial seed by a finite sequence of mutations, and each mutation passes between adjacent $g$-cones.
	\end{proof}
	
	\begin{proposition}\label{prop:gfin-gcomp}
		Let $Q$ satisfy Assumption~\ref{ass:gfan}. If $Q$ is $g$-finite, then it is $g$-complete.
	\end{proposition}
	
	\begin{proof}
		We follow the finite-fan argument in \cite[Proposition~4.8]{A21}. The assertion is immediate in rank $1$, so assume that $n\geq2$. Suppose that $Q$ is $g$-finite, and set
		\[
		F_1:=|\mathcal F(Q)|,
		\qquad
		F_2:=\bigcup_{g\in\gseed_{n-2}(Q)}C(g).
		\]
		First, we prove that $F_1 \setminus F_2$ is an open and closed set in $\mathbb{R}^n \setminus F_2$. Since there are only finitely many $g$-cones, $F_1 \setminus F_2$ is closed. Every cone of dimension at most $n-2$ is contained in one of dimension $n-2$, since it is generated by a subset of a $g$-cluster. Hence every point of $F_1\setminus F_2$ lies either in the interior of a $g$-cone or in the relative interior of a wall. In the second case, Lemma~\ref{lem:two-cones} supplies a neighborhood contained in $F_1$. Hence $F_1 \setminus F_2$ is open in $\mathbb R^n\setminus F_2$.
		
		Since $F_2$ is a finite union of cones of dimension at most $n-2$, $\mathbb R^{n}\setminus F_2$ is connected. Thus its nonempty open and closed subset $F_1\setminus F_2$ equals $\mathbb R^{n}\setminus F_2$, and hence $F_1=\mathbb R^n$.
	\end{proof}

\subsection{Mutation of the initial quiver}
  For $k\in Q_0$, let $H_k^{+} = \{\mathbf{x} = (x_1,\dots,x_n)^{\mathsf T} \mid x_k \geq 0\}$, $H_k^{-} = \{\mathbf{x} = (x_1,\dots,x_n)^{\mathsf T} \mid x_k \leq 0\}$, and $H_k = H_k^+ \cap H_k^{-}$. Let $\mu_k^Q:\mathbb R^{n}\to\mathbb R^{n}$, $(x_1,\dots,x_n)^{\mathsf T} \mapsto (x'_1,\dots,x'_n)^{\mathsf T}$ be the map defined by
  \begin{equation}\label{eq:initial-mutation}
    x'_i=
    \begin{cases}
      -x_k & \text{if $i=k$},\\
      x_i+[b_{ik}^Q]_+x_k+b_{ik}^Q[-x_k]_{+}
      & \text{if $i\neq k$}.
    \end{cases}
  \end{equation}

  \begin{lemma}\label{lem:homofg-fanmut}
    The map $\mu_k^Q$ is a piecewise-linear homeomorphism that is linear on each of $H_k^{+}$ and $H_k^{-}$. Moreover, its inverse is given by $\mu_k^{\mu_k(Q)}$.
  \end{lemma}
  \begin{proof}
    By direct calculation, we see that $\mu_k^Q$ is linear on each of the closed domains $H_k^{+}$ and $H_k^{-}$. Since these two domains intersect at the hyperplane $H_k$, $\mu_k^Q$ is continuous.

    Since $b_{ik}^{\mu_k(Q)}=-b_{ik}^Q$ for $i\neq k$, direct substitution in \eqref{eq:initial-mutation} gives
		\[
		\mu_k^{\mu_k(Q)}\circ\mu_k^Q=\mathrm{id}.
		\]
		Interchanging $Q$ and $\mu_k(Q)$ gives the reverse identity $\mu_k^Q \circ \mu_k^{\mu_k(Q)} = \mathrm{id}$.
  \end{proof}
  This map gives the formula of the mutation of initial quivers under Assumption~\ref{ass:gfan}.
  \begin{lemma}\label{lem:initial-seed-mut}
    Let $Q$ satisfy Assumption~\ref{ass:gfan}. Let $(\mathbf{g}_1,\ldots,\mathbf{g}_n)$ be the ordered tuple of $g$-vectors obtained from the initial seed for $Q$ by applying mutations along $k_1,\dots,k_r$. Then for each $k\in Q_0$, the ordered tuple of $g$-vectors obtained from the initial seed for $\mu_k(Q)$ by applying mutations along $k,k_1,\ldots,k_r$ is $\bigl(\mu_k^Q(\mathbf{g}_1),\ldots,\mu_k^Q(\mathbf{g}_n)\bigr)$. In particular, applying $\mu_k^Q$ to each vector gives a bijection between the $g$-clusters of $Q$ and those of $\mu_k(Q)$.
  \end{lemma}
  \begin{proof}
    The initial-seed mutation formula \cite[Proposition~6.1(c), equation~(6.2)]{AC25}, applied to each column, is precisely \eqref{eq:initial-mutation}.
  \end{proof}

  For a sequence $\mathbf{k} = (k_1,\dots,k_r)$ of vertices of $Q$, set $Q^0 = Q$ and $Q^{j} = \mu_{k_{j}}(Q^{j-1})$ for $1\leq j\leq r$. We define
  \[
  \Phi_{\mathbf{k}} = \mu_{k_r}^{Q^{r-1}} \circ \cdots \mu_{k_2}^{Q^1} \circ \mu_{k_1}^Q,
  \]
  with $\Phi_{\mathbf{k}}=\mathrm{id}$ when $r=0$.

  \begin{proposition}\label{prop:comp-of-initialmutation}
    Let $Q$ satisfy Assumption~\ref{ass:gfan}. Let $(\mathbf{g}_1,\ldots,\mathbf{g}_n)$ be the ordered tuple of $g$-vectors obtained from the initial seed for $Q$ by applying mutations along $\mathbf{k}=(k_1,\ldots,k_r)$. Then we have $\Phi_{\mathbf{k}}(\mathbf{g}_i) = \mathbf{e}_i$ for any $i\in Q_0$.
  \end{proposition}

  \begin{proof}
    Since each mutation is an involution, the assertion follows from Lemma~\ref{lem:initial-seed-mut}.
  \end{proof}

  For the $g$-fans, this initial seed mutation is viewed as follows.
	\begin{proposition}\label{prop:gfan-mut}
		Let $Q$ satisfy Assumption~\ref{ass:gfan}. The homeomorphism $\mu_k^Q$ is linear on every cone of $\mathcal F(Q)$ and sends these cones bijectively onto the cones of $\mathcal F(\mu_k(Q))$. In particular, $g$-finiteness, $g$-completeness, and $g$-density are mutation-invariant under Assumption~\ref{ass:gfan}.
	\end{proposition}
	
	\begin{proof}
		By Lemma~\ref{lem:assumption-restriction}, $\mu_k(Q)$ also satisfies Assumption~\ref{ass:gfan}. For any $g$-cluster $g = \{\mathbf{g}_1,\dots,\mathbf{g}_n\}$, by sign-coherence, the $k$th coordinates of the vectors in $g$ are either all nonnegative or all nonpositive. Hence $C(g)$ lies in one of $H_k^{+}$ and $H_k^{-}$. By Lemma~\ref{lem:initial-seed-mut}, $\{\mu_k^Q(\mathbf g_1),\ldots,\mu_k^Q(\mathbf g_n)\}$ is a $g$-cluster for $\mu_k(Q)$. Since $\mu_k^Q$ is linear on the half-space containing $C(g)$, it sends $C(g)$ to the cone generated by this cluster. The same applies to every face of $C(g)$.
		
		Since $\mu_k^Q$ and its inverse are continuous by Lemma~\ref{lem:homofg-fanmut}, $\mu_k^Q$ is a homeomorphism carrying $|\mathcal F(Q)|$ onto $|\mathcal F(\mu_k(Q))|$. Finiteness of the set of $g$-clusters is preserved by the bijection, and completeness and density of the support are preserved by the homeomorphism.
	\end{proof}
  
\subsection{Coordinate projections of \texorpdfstring{$g$}{g}-cones}
  From now on, we naturally identify $\mathbb{R}^{Q_0} = \mathbb{R}^n$ for simplicity. For a nonempty subset $I\subseteq Q_0$, let $\pi_I:\mathbb R^{Q_0}\to\mathbb R^I$ be the coordinate projection. By Lemma~\ref{lem:assumption-restriction}, $Q_I$ also satisfies Assumption~\ref{ass:gfan}.
	
	\begin{lemma}\label{lem:pi-commute}
		Let $Q$ satisfy Assumption~\ref{ass:gfan}. For $k\in I$, we have
		\[
		\pi_I\circ\mu_k^Q=\mu_k^{Q_I}\circ\pi_I.
		\]
	\end{lemma}
	
	\begin{proof}
		For $i\in I$, the $i$th component of \eqref{eq:initial-mutation} depends only on $x_i$, $x_k$, and $b_{ik}^Q$, all of which are unchanged by restriction to $I$.
	\end{proof}
	
	The image of a $g$-cone under $\pi_I$ is a closed polyhedral cone of full dimension in $\mathbb R^I$. Indeed, it is finitely generated by the projections of a $g$-cluster, so it is polyhedral and closed. These generators span $\mathbb R^I$, because the $g$-cluster is a basis by \eqref{eq:GC-duality}. We use this fact, and the fact that a full-dimensional closed convex cone is the closure of its interior, without further mention below.
	
	\begin{lemma}\label{lem:pi-noncross-wall}
		Let $Q$ satisfy Assumption~\ref{ass:gfan}. Let $C$ be a $g$-cone of $Q$, and $W$ a wall of $\mathcal F(Q_I)$. Then
		\[
		  \int(\pi_I(C))\cap W=\varnothing.
		\]
	\end{lemma}
	
	\begin{proof}
		Choose a seed for $Q_I$ whose $g$-cone contains $W$ as a wall. Let $(\mathbf g_i)_{i\in I}$ be its ordered tuple of $g$-vectors, and choose a sequence $k_1,\ldots,k_r$ of mutations from the initial seed for $Q_I$ to this seed. Set $Q^0=Q$, $Q^j=\mu_{k_j}(Q^{j-1})$ for $1\leq j\leq r$, and $R=Q^r$. By \eqref{eq:restriction-mutation}, applying the same sequence of mutations to $Q_I$ gives $R_I$. Write the compositions of piecewise linear homeomorphisms in Lemma~\ref{lem:homofg-fanmut} as $\Phi = \mu_{k_{r}}^{Q^{r-1}} \circ \cdots \circ \mu_{k_1}^{Q^0}$ and $\Phi_I = \mu_{k_{r}}^{Q_{I}^{r-1}} \circ \cdots \circ \mu_{k_1}^{Q_{I}^0}$. By Proposition~\ref{prop:gfan-mut}, applied to each cone, these maps induce bijections $\Phi: \mathcal{F}(Q) \to \mathcal{F}(R)$ and $\Phi_I:\mathcal{F}(Q_I) \to \mathcal{F}(R_I)$. By Proposition~\ref{prop:comp-of-initialmutation}, we have $\Phi_I(\mathbf{g}_i) = \mathbf{e}_i$ for any $i \in I$. In particular, $\Phi_I(W)$ is contained in $H_{k;I} := \{\mathbf{x} \in \mathbb{R}^I \mid x_k = 0 \}$ for some $k \in I$. By Lemma~\ref{lem:pi-commute}, applied at each intermediate quiver,
		\[
		\pi_I\circ\Phi=\Phi_I\circ\pi_I.
		\]
		
		If $x\in\int(\pi_I(C))\cap W$, then the homeomorphism $\Phi_I$ gives
		\[
		\Phi_I(x)\in\int\bigl(\Phi_I(\pi_I(C))\bigr)\cap H_{k;I}
		=\int\bigl(\pi_I(\Phi(C))\bigr)\cap H_{k;I}.
		\]
		On the other hand, since $\Phi(C)$ is a $g$-cone of $\mathcal{F}(R)$, by sign-coherence, its projection lies in one of the closed half-spaces bounded by $H_{k;I}$. The interior of such a set cannot meet $H_{k;I}$. This is a contradiction.
	\end{proof}
	
	The next proof uses an elementary convexity observation. Let $P\subseteq\mathbb R^I$ be a full-dimensional closed convex cone such that
	\[
	\int(P)\cap W=\varnothing
	\quad\text{for every wall }W\text{ of }\mathcal F(Q_I).
	\]
	If $\int(P)$ meets the interior of a $g$-cone $\sigma$, then $P\subseteq\sigma$. Otherwise, since $P=\overline{\int(P)}$ and $\sigma$ is closed, there would be a point of $\int(P)$ outside $\sigma$. The line segment joining it to a point of $\int(P)\cap\int(\sigma)$ stays in $\int(P)$ and meets the boundary $\partial\sigma$. The boundary is a union of walls, giving a contradiction. This argument uses only the finitely many facets of $\sigma$.
	
	\begin{lemma}\label{lem:pi-adjacent-cones}
		Let $Q$ satisfy Assumption~\ref{ass:gfan}. Let $C$ and $C'$ be adjacent $g$-cones of $Q$. If $\pi_I(C)\subseteq\sigma$ for a $g$-cone $\sigma$ of $Q_I$, then $\pi_I(C')$ is contained in $\sigma$ or in a $g$-cone adjacent to $\sigma$.
	\end{lemma}
	
	\begin{proof}
		Let $F=C\cap C'$ be the common wall. Put $r = |I|$. Since $\dim\ker\pi_I=n-r$, we have
		\[
		\dim\pi_I(F)\ge (n-1)-(n-r)=r-1.
		\]
		Also, $\pi_I(F)\subseteq \sigma\cap\pi_I(C')$. If $\pi_I(F)$ meets $\operatorname{int}(\sigma)$, then so does $\operatorname{int}(\pi_I(C'))$, since $\pi_I(C')$ is the closure of its interior. Lemma~\ref{lem:pi-noncross-wall} and the convexity observation above give $\pi_I(C')\subseteq\sigma$.
		
		Suppose that $\pi_I(F)\subseteq\partial\sigma$. Since $\pi_I(F)$ is a convex cone of dimension at least $r-1$, it cannot be contained in the union of the faces of $\sigma$ of dimension at most $r-2$. Thus, we can choose a point $x\in\pi_I(F)$ lying in the relative interior of a wall $\tau$ of $\sigma$. Let $\sigma'$ be the other $g$-cone containing $\tau$. By Lemma~\ref{lem:two-cones}, some neighborhood $U$ of $x$ is contained in $\sigma\cup\sigma'$. Since $x\in\pi_I(C')$ and $\pi_I(C')$ is the closure of its interior, we can choose
		\[
		y\in U\cap\operatorname{int}(\pi_I(C')).
		\]
		Since $y$ lies on no wall of $\mathcal F(Q_I)$ by Lemma~\ref{lem:pi-noncross-wall}, it lies in $\operatorname{int}(\sigma)$ or $\operatorname{int}(\sigma')$. The convexity observation above gives $\pi_I(C')\subseteq\sigma$ or $\pi_I(C')\subseteq\sigma'$, as required.
	\end{proof}
	
	\begin{proposition}\label{prop:pi-in-cone}
		Let $Q$ satisfy Assumption~\ref{ass:gfan}. For every $g$-cone $C$ of $Q$, there is a $g$-cone $\sigma$ of $Q_I$ such that $\pi_I(C)\subseteq\sigma$.
	\end{proposition}
	
	\begin{proof}
		By Lemma~\ref{lem:two-cones}, there is a sequence of adjacent $g$-cones
		\[
		C_0,C_1,\ldots,C_m=C
		\]
		starting at the initial cone $C_0=\mathbb R_{\geq0}^{Q_0}$. Since its projection is the initial cone $\mathbb R_{\geq0}^I$ of $Q_I$, the assertion follows by applying Lemma~\ref{lem:pi-adjacent-cones} successively along the sequence.
	\end{proof}
	
	\begin{theorem}\label{thm:I-gcomp-gdense}
		Let $Q$ satisfy Assumption~\ref{ass:gfan}, and $I\subseteq Q_0$ be nonempty. Then
		\begin{equation}\label{eq:projected-support}
			\pi_I\bigl(|\mathcal F(Q)|\bigr)=|\mathcal F(Q_I)|.
		\end{equation}
		In particular, the following hold.
		\begin{enumerate}[label=\textup{(\arabic*)}]
			\item If $Q$ is $g$-complete, then $Q_I$ is $g$-complete.
			\item If $Q$ is $g$-dense, then $Q_I$ is $g$-dense.
		\end{enumerate}
	\end{theorem}
	
	\begin{proof}
		Proposition~\ref{prop:pi-in-cone} gives the inclusion $\pi_I\bigl(|\mathcal F(Q)|\bigr)\subseteq|\mathcal F(Q_I)|$. For the reverse inclusion, let $\sigma$ be a $g$-cone of $Q_I$. Choose a mutation sequence from the initial seed for $Q_I$ to a seed with $g$-cone $\sigma$, and apply the same sequence to the initial seed for $Q$. By the proof of Lemma~\ref{lem:SC-full}, the resulting $g$-vectors indexed by $I$ project to the generators of $\sigma$. For $j\notin I$, the $j$th $g$-vector is $\mathbf{e}_j$ and projects to zero. Thus the resulting $g$-cone of $Q$ projects onto $\sigma$, proving \eqref{eq:projected-support}.
		
		If $Q$ is $g$-complete, then \eqref{eq:projected-support} gives $|\mathcal F(Q_I)|=\pi_I(\mathbb R^{Q_0})=\mathbb R^I$. If $Q$ is $g$-dense, then for each nonempty open set $U\subseteq\mathbb R^I$, the set $\pi_I^{-1}(U)$ is nonempty and open, and hence it meets $|\mathcal F(Q)|$. By \eqref{eq:projected-support}, $U$ meets $|\mathcal F(Q_I)|$.
	\end{proof}

\subsection{Rank-\texorpdfstring{$2$}{2} obstructions and completeness}	
	\begin{lemma}\label{lem:rank2-not-gcomp}
		Let $Q$ be a weighted quiver of rank $2$ whose arrow has weight at least $2$. Then $Q$ satisfies Assumption~\ref{ass:gfan}, but is not $g$-complete.
	\end{lemma}
	
	\begin{proof}
		Since $Q$ has the totally DB2-property, Theorem~\ref{thm:SC-D-rank2} and \cite[Theorem~7.7]{AC25} give Assumption~\ref{ass:gfan}. By \cite[Example~7.8]{AC25} (see also \cite[Proposition~9.6]{Re14}), $Q$ is not $g$-complete.
	\end{proof}
	
	\begin{proposition}\label{prop:gcomp-2fin}
		Let $Q$ satisfy Assumption~\ref{ass:gfan}. If $Q$ is $g$-complete, then $Q$ is $2$-finite.
	\end{proposition}
	
	\begin{proof}
		If $Q$ is not $2$-finite, then there is a quiver $Q'$ in the mutation class of $Q$ such that it has an arrow of weight at least $2$. By Theorem~\ref{thm:I-gcomp-gdense}(1) and Lemma~\ref{lem:rank2-not-gcomp}, $Q'$ is not $g$-complete, and hence $Q$ is also not $g$-complete by Proposition~\ref{prop:gfan-mut}.
	\end{proof}
	
	\begin{theorem}\label{thm:main-classification}
		Let $Q$ satisfy Assumption~\ref{ass:gfan}. Then the following conditions are equivalent:
		\begin{enumerate}[label=\textup{(\arabic*)}]
			\item $Q$ is $g$-finite;
			\item $Q$ is $g$-complete;
			\item $Q$ is $2$-finite;
			\item $Q$ is of finite Coxeter type.
		\end{enumerate}
	\end{theorem}
	
	\begin{proof}
		The assertion follows from Theorem~\ref{thm:gfin-Cox}, Propositions~\ref{prop:gfin-gcomp} and~\ref{prop:gcomp-2fin}.
	\end{proof}

\section{Proof of Theorem~\ref{thm:2fin-Cox}}\label{sec:pf-2fin-Cox}
	Let $Q$ be a weighted quiver. To prove Theorem~\ref{thm:2fin-Cox}, we follow the method of \cite[Section~9]{FZ03a}: we first treat trees and chordless cycles, and then use induction on the number of vertices. Throughout this section, we use without further mention the basic facts about $2$-finiteness and the totally DB2-property stated in Section~\ref{sec:2fin-DB2}. We will use the identity $\bra{5}^2=\bra{5}+1$. In mutation diagrams,
	\begin{tikzpicture}
		\node[myMQV] at (0,0) {};
	\end{tikzpicture}
	indicates the vertex at which we mutate. Throughout this section, a weight of zero means that the corresponding arrow is absent.

\subsection{Trees, cycles, and reorientations}
	A \emph{chordless cycle} in a quiver is a full subquiver whose underlying graph is a cycle. Mutation at a source or sink only reverses the incident arrows. We record a reorientation lemma and introduce two families of quivers used below.
	
	\begin{lemma}[cf.~{\cite[Proposition~9.2]{FZ03a}}]\label{lem: chain reorientations}
		Let $T$ be a full subtree of the underlying weighted graph of $Q$. Suppose that $T$ meets the rest of the graph only at a vertex $v$: no vertex of $T\setminus\{v\}$ is adjacent to a vertex outside $T$. Then any reorientation of the arrows of $Q$ whose endpoints lie in $T$, with all other arrows unchanged, can be obtained by source and sink mutations at vertices of $T\setminus\{v\}$. In particular, it is mutation-equivalent to $Q$.
	\end{lemma}
	
	\begin{proof}
		Since the proof of \cite[Proposition~9.2]{FZ03a} uses only source and sink mutations, it applies unchanged to arbitrary positive real weights.
	\end{proof}
	
	We call the operation in Lemma~\ref{lem: chain reorientations} a \emph{tree reorientation}. In particular, all orientations of a weighted tree are mutation-equivalent.
	
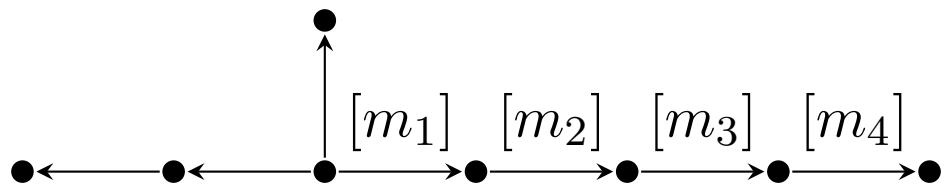
\begin{figure}[htbp]
		\centering
		\begin{tikzpicture}
			\node[myQV] (0) at (0,0) {};
			\node[myQV] (-1) at (-1,0) {};
			\node[myQV] (-2) at (-2,0) {};
			\node[myQV] (A) at (0,1) {};
			\draw[->] (0) -> (-1);
			\draw[->] (-1) -> (-2);
			\draw[->] (0) -> (A);
			\foreach \x [evaluate = \x as \y using int(\x - 1)] in {1,2,3,4}
			{
				\node[myQV] (\x) at (\x,0) {};
				\draw[->] (\y) -> node [pos = 0.5, above] {$\bra{m_{\x}}$} (\x);
			}
		\end{tikzpicture}
		\caption{Tree quiver $T_{4,1,2}(m_1,m_2,m_3,m_4)$.}
		\label{fig: tree quiver}
	\end{figure}

	\begin{definition}
		Let $a_1,a_2,a_3\in\mathbb Z_{\ge0}$ and $\mathbf m=(m_1,\ldots,m_{a_1})\in\mathbb Z_{\ge3}^{a_1}$. The weighted quiver $T_{a_1,a_2,a_3}(\mathbf m)$ has a central vertex $v$ and three directed paths
		\[
		v=v_0^i\longrightarrow v_1^i\longrightarrow\cdots
		\longrightarrow v_{a_i}^i
		\qquad(i=1,2,3),
		\]
		which have no other vertices in common. The $k$th arrow on the first path has weight $\bra{m_k}$, and all other arrows have weight $1$. We omit $\mathbf m$ when all weights are $1$. See Figure~\ref{fig: tree quiver} for an example.
	\end{definition}
	
\begin{figure}[htbp]
		\centering
		\begin{tikzpicture}
			\node[myQV] (0) at (-30:2) {};
			\node[myQV] (1) at (0:2) {};
			\node[myQV] (2) at (30:2) {};
			\node at (0) [above left] {$v_1$};
			\node at (1) [left] {$v_2$};
			\node at (2) [below left] {$v_3$};
			\foreach \x in {3,...,11}
			{
				\node[myQV] (\x) at (30*\x - 30: 2) {};
			}
			\foreach \x [evaluate = \x as \y using int(\x+1)] in {0,...,10}
			{
				\draw[->] (\x) -> (\y);
			}
			\draw[->] (11) -> (0);
			\node[myQV] (A1) at ($(0)+(-30: 0.9)$) {};
			\node[myQV] (A2) at ($(0)+(-30: 0.9*2)$) {};
			\node[myQV] (B1) at ($(1)+(0: 0.9)$) {};
			\node[myQV] (B2) at ($(1)+(0: 0.9*2)$) {};
			\node[myQV] (B3) at ($(1)+(0: 0.9*3)$) {};
			\node[myQV] (C1) at ($(2)+(30: 0.9)$) {};
			
			\draw[->] (A1) -> node [pos = 0.5, above right=-3pt] {\tiny $\bra{m_1}$} (0);
			\draw[->] (A2) -> node [pos = 0.5, above right=-3pt] {\tiny $\bra{m_2}$} (A1);
			\draw[->] (B1) -> (1);
			\draw[->] (B2) -> (B1);
			\draw[->] (B3) -> (B2);
			\draw[->] (C1) -> (2);
			
			\draw [
			thick,
			decoration={
				brace,
				mirror,
				raise=3pt
			},
			decorate
			] (0.south east) -- (A2.south east) node [pos=0.5, anchor=north east, xshift = 5pt, yshift = -7pt] {$a_1$ vertices};
			\draw [
			thick,
			decoration={
				brace,
				mirror,
				raise=3pt
			},
			decorate
			] (1.south east) -- (B3.south east) node [pos=0.5, anchor=north, yshift = -5pt] {$a_2$ vertices};
			\draw [
			thick,
			decorate,
			decoration={
				brace,
				mirror,
				raise=3pt
			},
			decorate
			] (2.south east) -- (C1.south east) node [pos=0.5, anchor = north west, yshift = -2pt] {$a_3$ vertices};
		\end{tikzpicture}
		\caption{The quiver $S_{3,4,2}^{9}(m_1,m_2)$.}
		\label{fig: example of S-quiver}
	\end{figure}
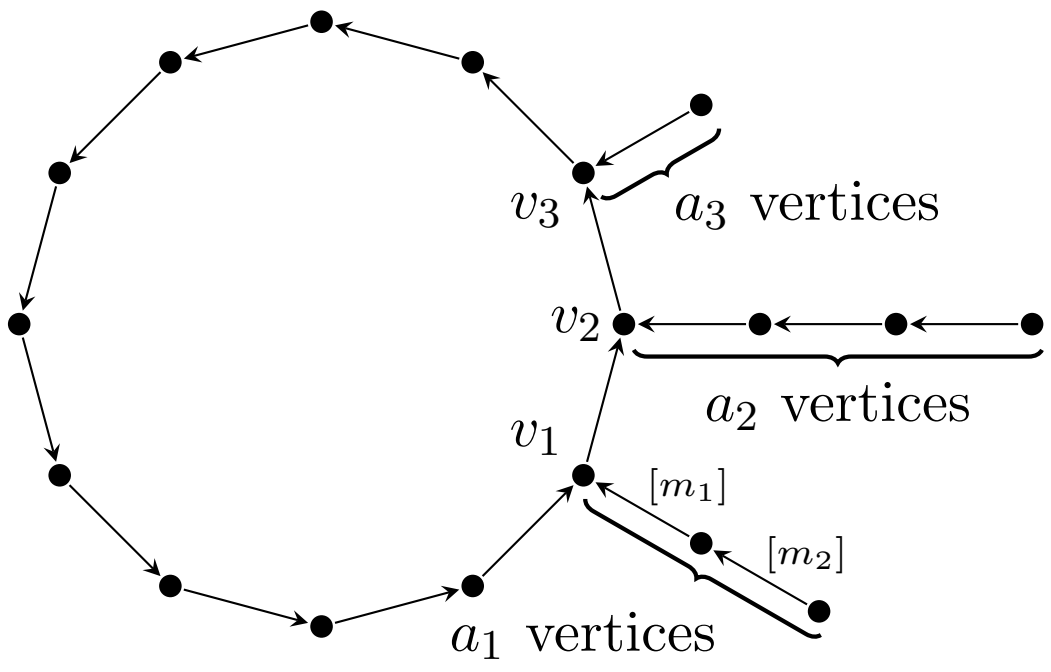

	\begin{definition}\label{def:quiverS}
		Let $a_1,a_2,a_3\in\mathbb Z_{>0}$, $s\in\mathbb Z_{\ge0}$, and $\mathbf m=(m_1,\ldots,m_{a_1-1}) \in\mathbb Z_{\ge3}^{a_1-1}$. The weighted quiver $S_{a_1,a_2,a_3}^{s}(\mathbf m)$ is obtained from the directed cycle
		\[
		v_1\longrightarrow v_2\longrightarrow\cdots
		\longrightarrow v_{s+3}\longrightarrow v_1
		\]
		by attaching a directed path
		\[
		u_{a_i}^i\longrightarrow\cdots\longrightarrow u_2^i
		\longrightarrow u_1^i=v_i
		\]
		at each $v_i$ for $i=1,2,3$. The arrow $u_{k+1}^1\to u_k^1$ has weight $\bra{m_k}$, and all other arrows have weight $1$. We omit $\mathbf m$ when all weights are $1$. A path with $a_i=1$ adds no vertex to the cycle. See Figure~\ref{fig: example of S-quiver} for an example.
	\end{definition}
	
	\begin{lemma}[cf.~{\cite[Lemma~9.6]{FZ03a}}]\label{lem: mutation equivalence between S and T}

		With the notation above, $S_{a_1,a_2,a_3}^{s}(\mathbf m)$ is mutation-equivalent to
		\[
		T_{a_1+a_3-1,a_2,s}(\widetilde{\mathbf m}),
		\qquad
		\widetilde{\mathbf m}=(\underbrace{3,\ldots,3}_{a_3\textup{ entries}},\mathbf m).
		\]
	\end{lemma}
	
	\begin{proof}
    We label the vertices of $S^{s}_{a_1,a_2,a_3}(\mathbf{m})$ as in Definition~\ref{def:quiverS}.
		Mutating it successively at
		\[
		v_{s+3},v_{s+2},\ldots,v_3,u^3_2,\ldots,u^3_{a_3}
		\]
		gives a reorientation of $T_{a_1+a_3-1,a_2,s}(\widetilde{\mathbf m})$. The assertion follows from Lemma~\ref{lem: chain reorientations}.
	\end{proof}

\subsection{Quivers of rank \texorpdfstring{$3$}{3}}
	We first treat quivers of rank $3$ and prove the following classification.
	
	\begin{proposition}\label{prop: 2-finte rank 3 quivers}
		A connected quiver $Q\in\D$ of rank $3$ is $2$-finite if and only if, up to relabeling and tree reorientation, it is one of the quivers in Table~\ref{tab: 2-finite of rank 3}. 
	\end{proposition}
	
	\begin{table}[htbp]
		\centering
		\renewcommand{\arraystretch}{1.25}
		\begin{tabular}{c|ccc}
			\rule[-20pt]{0pt}{50pt}
			type $A_3$
			&
			\begin{tikzpicture}[baseline = {(0,0)}]
				\node[myQV] (1) at (1,0) {};
				\node[myQV] (2) at (2,0) {};
				\node[myQV] (3) at (3,0) {};
				\draw[->] (1) -> (2);
				\draw[->] (2) -> (3);
			\end{tikzpicture}
			&
			\begin{tikzpicture}[baseline = {(0,0)}, scale = 0.75]
				\node[myQV] (1) at (210: 1) {};
				\node[myQV] (2) at (330: 1) {};
				\node[myQV] (3) at (90: 1) {};
				\draw[->] (1) -> (2);
				\draw[->] (2) -> (3);
				\draw[->] (3) -> (1);
			\end{tikzpicture}
			&
			\\
			\hline
			\rule[-20pt]{0pt}{50pt}
			type $B_3$
			&
			\begin{tikzpicture}[baseline = {(0,0)}]
				\node[myQV] (1) at (1,0) {};
				\node[myQV] (2) at (2,0) {};
				\node[myQV] (3) at (3,0) {};
				\draw[->] (1) -> (2);
				\draw[->] (2) -> (3) node [pos = 0.5, above] {$\bra{4}$};
			\end{tikzpicture}
			&
			\begin{tikzpicture}[baseline = {(0,0)}, scale = 0.75]
				\node[myQV] (1) at (210: 1) {};
				\node[myQV] (2) at (330: 1) {};
				\node[myQV] (3) at (90: 1) {};
				\draw[->] (1) -> (2);
				\draw[->] (2) -> (3) node [pos = 0.5, right] {$\bra{4}$};
				\draw[->] (3) -> (1) node [pos = 0.5, left] {$\bra{4}$};
			\end{tikzpicture}
			&
			\\
			\hline
			\rule[-20pt]{0pt}{50pt}
			type $H_3$
			&
			\begin{tikzpicture}[baseline = {(0,0)}]
				\node[myQV] (1) at (1,0) {};
				\node[myQV] (2) at (2,0) {};
				\node[myQV] (3) at (3,0) {};
				\draw[->] (1) -> (2);
				\draw[->] (2) -> (3) node [pos = 0.5, above] {$\bra{5}$};
			\end{tikzpicture}
			&
			\begin{tikzpicture}[baseline = {(0,0)}, scale = 0.75]
				\node[myQV] (1) at (210: 1) {};
				\node[myQV] (2) at (330: 1) {};
				\node[myQV] (3) at (90: 1) {};
				\draw[->] (1) -> (2);
				\draw[->] (2) -> (3) node [pos = 0.5, right] {$\bra{5}$};
				\draw[->] (3) -> (1) node [pos = 0.5, left] {$\bra{5}$};
			\end{tikzpicture}
			&
			\begin{tikzpicture}[baseline = {(0,0)}, scale = 0.75]
				\node[myQV] (1) at (210: 1) {};
				\node[myQV] (2) at (330: 1) {};
				\node[myQV] (3) at (90: 1) {};
				\draw[->] (1) -> (2) node [pos = 0.5, below] {$\bra{5}$};
				\draw[->] (2) -> (3) node [pos = 0.5, right] {$\bra{5}$};
				\draw[->] (3) -> (1) node [pos = 0.5, left] {$\bra{5}$};
			\end{tikzpicture}
		\end{tabular}
		\caption{The $2$-finite quivers of rank $3$ in $\D$.}
		\label{tab: 2-finite of rank 3}
	\end{table}
	
	To prove Proposition~\ref{prop: 2-finte rank 3 quivers}, we prepare three lemmas.
	
	\begin{lemma}\label{lem: rank 3 tree case}
		If a connected quiver $Q\in\D$ of rank $3$ is acyclic and $2$-finite, then it is a finite Coxeter quiver of type $A_3$, $B_3$, or $H_3$. In particular, it is a tree with weights at most $\bra{5}$.
	\end{lemma}
	
	\begin{proof}
		All weights of $Q$ are at least $1$. If its underlying graph is a triangle, then the following mutation gives the weight $ab+c\geq2$, a contradiction:
		
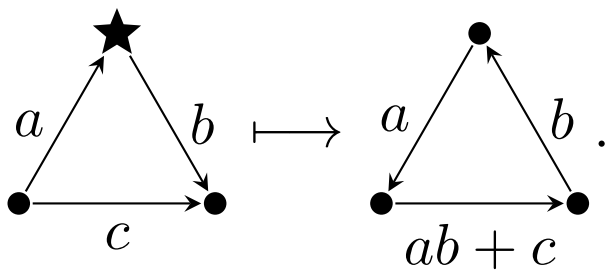
\[
		\begin{tikzpicture}[baseline = {(0,0)}, scale = 0.75]
			\node[myMQV] (1) at (90: 1) {};
			\node[myQV] (2) at (210: 1) {};
			\node[myQV] (3) at (330: 1) {};
			\draw[->] (2) -> (1) node [pos = 0.5, left] {$a$};
			\draw[->] (2) -> (3) node [pos = 0.5, below] {$c$};
			\draw[->] (1) -> (3) node [pos = 0.5, right] {$b$};
		\end{tikzpicture}
		\longmapsto
		\begin{tikzpicture}[baseline = {(0,0)}, scale = 0.75]
			\node[myQV] (1) at (90: 1) {};
			\node[myQV] (2) at (210: 1) {};
			\node[myQV] (3) at (330: 1) {};
			\draw[->] (1) -> (2) node [pos = 0.5, left] {$a$};
			\draw[->] (2) -> (3) node [pos = 0.5, below] {$ab+c$};
			\draw[->] (3) -> (1) node [pos = 0.5, right] {$b$};
		\end{tikzpicture}.
		\]
		Thus $Q$ is a tree. By Lemma~\ref{lem: chain reorientations}, we may reorient it as a directed path with weights $a\geq b$. Mutation at the middle vertex gives the weight $ab$. Since $b^2\leq ab<2$, we have $b=1$. Then the following two mutations give the weight $a^2-1$:
		\[
		\begin{tikzpicture}[baseline = {(0,0)}, scale = 0.75]
			\node[myMQV] (1) at (90: 1) {};
			\node[myQV] (2) at (210: 1) {};
			\node[myQV] (3) at (330: 1) {};
			\draw[->] (2) -> (1) node [pos = 0.5, left] {$a$};
			\draw[->] (1) -> (3);
		\end{tikzpicture}
		\longmapsto
		\begin{tikzpicture}[baseline = {(0,0)}, scale = 0.75]
			\node[myQV] (1) at (90: 1) {};
			\node[myMQV] (2) at (210: 1) {};
			\node[myQV] (3) at (330: 1) {};
			\draw[->] (1) -> (2) node [pos = 0.5, left] {$a$};
			\draw[->] (2) -> (3) node [pos = 0.5, below] {$a$};
			\draw[->] (3) -> (1);
		\end{tikzpicture}
		\longmapsto
		\begin{tikzpicture}[baseline = {(0,0)}, scale = 0.75]
			\node[myQV] (1) at (90: 1) {};
			\node[myQV] (2) at (210: 1) {};
			\node[myQV] (3) at (330: 1) {};
			\draw[->] (2) -> (1) node [pos = 0.5, left] {$a$};
			\draw[->] (3) -> (2) node [pos = 0.5, below] {$a$};
			\draw[->] (1) -> (3) node [pos = 0.5, right] {$a^2 - 1$};
		\end{tikzpicture}.
		\]
		Since $a^2-1<2$, we obtain $a\in\{1,[4],[5]\}$.
	\end{proof}
	
	Next, we consider the case of directed triangles. By the obvious symmetry, we may assume that a directed triangle with weights in $\mathcal W_{<2}$ has the form
	\begin{equation}\label{eq: rank 3 cyclic quiver}
		Q:\quad
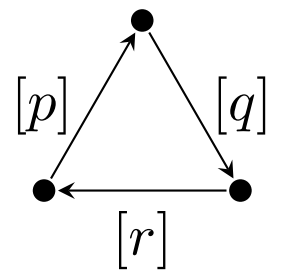
		\begin{tikzpicture}[baseline={(current bounding box.center)}, scale=0.75]
			\node[myQV] (1) at (90:1) {};
			\node[myQV] (2) at (210:1) {};
			\node[myQV] (3) at (330:1) {};
			\draw[->] (2) -- node[midway,left] {$\bra{p}$} (1);
			\draw[->] (1) -- node[midway,right] {$\bra{q}$} (3);
			\draw[->] (3) -- node[midway,below] {$\bra{r}$} (2);
		\end{tikzpicture},
		\qquad p\geq q\geq r\geq3,
	\end{equation}
	where $p,q,r$ are integers. Let $Q_\ell$ be obtained from $Q$ by $\ell$ alternating mutations at the two lower vertices, starting at the lower right vertex, with $Q_0=Q$. For $\ell\geq-1$, define
	\[
	f_\ell=\frac{\bra{q}\sin\bigl((\ell+1)\pi/r\bigr)-\bra{p}\sin(\ell\pi/r)}{\sin(\pi/r)}.
	\]
	
	\begin{lemma}[cf.~{\cite[proof of Lemma~3.3]{FT19}}]\label{lem: 2 direction mutations}
		Let $\ell\in\mathbb Z_{\geq1}$. If $Q_0,\ldots,Q_{\ell-1}$ are directed triangles, then the three weights of $Q_\ell$ are $\bra{r}$, $|f_{\ell-1}|$, and $|f_\ell|$. Moreover, $Q_\ell$ is a directed triangle if and only if $f_\ell>0$.
	\end{lemma}
	
	\begin{proof}
		We have $f_{-1}=\bra{p}$, $f_0=\bra{q}$, and $f_{\ell+1}=\bra{r}f_\ell-f_{\ell-1}$. While the triangle is directed, this is exactly the signed weight change under mutation. A positive new signed weight gives a directed triangle, a zero weight gives a tree, and a negative weight gives an acyclic triangle. Induction proves the assertion.
	\end{proof}
	
	\begin{lemma}\label{lem: existence of acyclic quivers in two directed mutations}
		At least one of $Q_1,\ldots,Q_{r-1}$ is acyclic.
	\end{lemma}
	
	\begin{proof}
		If none of $Q_1,\ldots,Q_{r-2}$ is acyclic, then Lemma~\ref{lem: 2 direction mutations} applies at step $r-1$ and shows that $Q_{r-1}$ is acyclic since $f_{r-1}=-\bra{p}<0$.
	\end{proof}
	
	\begin{proof}[Proof of Proposition~\ref{prop: 2-finte rank 3 quivers}]
		Suppose first that $Q$ is $2$-finite. The acyclic case follows from Lemma~\ref{lem: rank 3 tree case}. Otherwise, we may assume that $Q$ has the form \eqref{eq: rank 3 cyclic quiver} and use the notation above. Lemmas~\ref{lem: 2 direction mutations} and~\ref{lem: existence of acyclic quivers in two directed mutations} give an acyclic quiver in its mutation class whose weights still include $\bra{r}$. Therefore, Lemma~\ref{lem: rank 3 tree case} gives $r\leq5$.
		
		If $r=3$, then $f_1=\bra{q}-\bra{p}\leq0$. A negative value would give an acyclic triangle, contrary to Lemma~\ref{lem: rank 3 tree case}. Hence $p=q$, and $Q_1$ is a tree with weights $1$ and $\bra{p}$. Lemma~\ref{lem: rank 3 tree case} gives $p=q\in\{3,4,5\}$.
		
		If $r=4$, then $f_1=\bra{4}\bra{q}-\bra{p}>0$, since $\bra{q}\geq\bra{4}=\sqrt2$ and $\bra{p}<2$. However, $f_2=\bra{q}-\bra{4}\bra{p}<0$, and hence $Q_2$ is an acyclic triangle, contrary to Lemma~\ref{lem: rank 3 tree case}.
		
		If $r=5$, then $f_1=\bra{5}\bra{q}-\bra{p}>0$ and $f_2=\bra{5}(\bra{q}-\bra{p})\leq0$. Since a negative value would again contradict Lemma~\ref{lem: rank 3 tree case}, $p=q$ and $Q_2$ is a tree with weights $\bra{5}$ and $(\bra{5}-1)\bra{p}$. By Lemma~\ref{lem: rank 3 tree case}, $p=q=5$.
		
		Conversely, for $p\in\{3,4,5\}$, a directed triangle with weights $(1,\bra{p},\bra{p})$ mutates to a tree with weights $1$ and $\bra{p}$. A triangle with all weights $\bra{5}$ first mutates to one with weights $(1,\bra{5},\bra{5})$ since $\bra{5}^2-\bra{5}=1$, and then it mutates to a tree with weights $1$ and $\bra{5}$. These trees are the finite Coxeter quivers of rank $3$, which are $2$-finite by \eqref{eq:Cox-2fin-DB2}.
	\end{proof}

\subsection{Classification of trees}
	We next treat tree quivers and prove the following classification.
	
	\begin{proposition}\label{prop: tree Coxeter quivers and 2-finite}
		A tree quiver $Q\in\D$ is $2$-finite if and only if it is a finite Coxeter quiver.
	\end{proposition}
	
	To prove Proposition~\ref{prop: tree Coxeter quivers and 2-finite}, we first identify several $2$-infinite trees. Let $m,m'\in\mathbb Z_{\geq3}$. We consider the following tree quivers:
	\begin{align*}
		Y_n(m) &=\
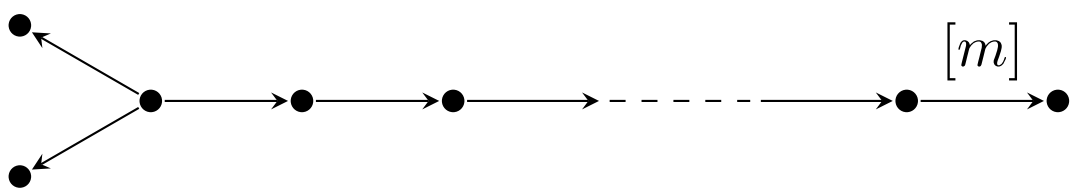
		\begin{tikzpicture}[baseline=-3pt]
			\node[myQV] (O) at (0,0) {};
			\node[myQV] (1) at (150:1) {};
			\node[myQV] (2) at (-150:1) {};
			\node[myQV] (3) at (1,0) {};
			\node[myQV] (4) at (2,0) {};
			\node[myCN] (5) at (3,0) {};
			\node[myCN] (6) at (4,0) {};
			\node[myQV] (7) at (5,0) {};
			\node[myQV] (8) at (6,0) {};
			\draw[->] (O) -> (1);
			\draw[->] (O) -> (2);
			\draw[->] (O) -> (3);
			\draw[->] (3) -> (4);
			\draw[->] (4) -> (5);
			\draw[dashed] (5) -- (6);
			\draw[->] (6) -> (7);
			\draw[->] (7) -> node [pos = 0.5, auto = left] {$\bra{m}$} (8);
		\end{tikzpicture}
		\quad
		(n \geq 4),
		\\
		\mathcal H_4 &=\
		
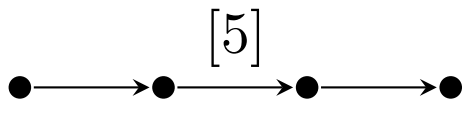
\begin{tikzpicture}[baseline=-3pt, scale = 0.95]
			\foreach \x in {1,2,3,4}
			{
				\node[myQV] (\x) at (\x,0) {};
			}
			\draw[->] (1) -> (2);
			\draw[->] (2) -> node [pos = 0.5, above] {$\bra{5}$} (3);
			\draw[->] (3) -> (4);
		\end{tikzpicture}\ ,
		\\
		Z_n(m,m') &=\
		
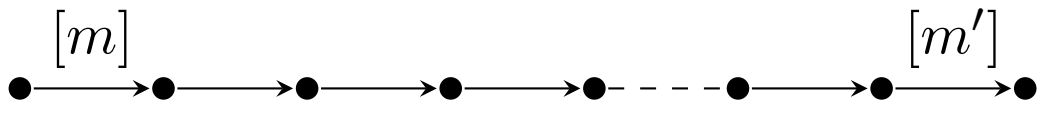
\begin{tikzpicture}[baseline=-3pt, scale = 0.95]
			\foreach \x in {1,...,8}
			{
				\node[myQV] (\x) at (\x,0) {};
			}
			\draw[->] (1) -> (2) node [pos = 0.5, above] {$\bra{m}$};
			\draw[->] (7) -> (8) node [pos = 0.5, above] {$\bra{m'}$};
			\foreach \x [evaluate = \x as \y using int(\x+1)] in {2,3,4,6}
			{
				\draw[->] (\x) -> (\y);
			}
			\draw[dashed] (5) -- (6);
		\end{tikzpicture}
		\quad
		(n \geq 3).
	\end{align*}
	By Lemma~\ref{lem: chain reorientations}, their $2$-finiteness does not depend on the orientation.
	
	\begin{lemma}\label{lem: special 2 infinite quivers}
		The quivers $Y_4(m)$ with $m\in\mathbb Z_{\geq4}$, $\mathcal H_4$, and $Z_5(3,5)$ are $2$-infinite.
	\end{lemma}
	
	\begin{proof}
		Reorient $Y_4(m)$ as on the left below. The following mutations give the weight $\bra{m}^2\geq 2$.
		
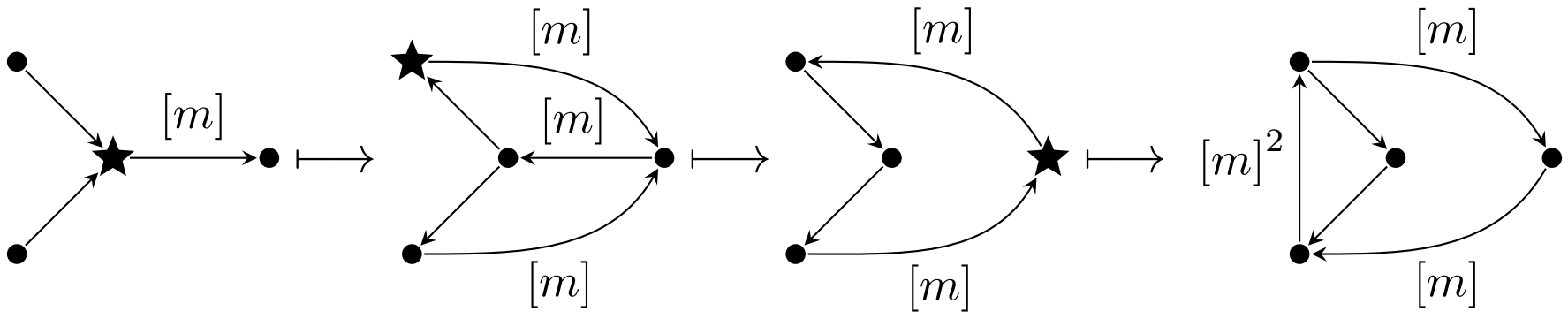
\[
		\begin{tikzpicture}[baseline=-3pt,scale=0.9,font=\small]
			\node[myQV] (q1) at (-0.8,0.8) {};
			\node[myQV] (q2) at (-0.8,-0.8) {};
			\node[myMQV] (q3) at (0,0) {};
			\node[myQV] (q4) at (1.3,0) {};
			\draw[->] (q1) -- (q3);
			\draw[->] (q2) -- (q3);
			\draw[->] (q3) -- node[midway,above] {$\bra{m}$} (q4);
		\end{tikzpicture}
		\longmapsto
		\begin{tikzpicture}[baseline=-3pt,scale=0.9,font=\small]
			\node[myMQV] (q1) at (-0.8,0.8) {};
			\node[myQV] (q2) at (-0.8,-0.8) {};
			\node[myQV] (q3) at (0,0) {};
			\node[myQV] (q4) at (1.3,0) {};
			\draw[->] (q3) -- (q1);
			\draw[->] (q3) -- (q2);
			\draw[->] (q4) -- node[above=-1pt,pos=0.6] {$\bra{m}$} (q3);
			\draw[->] (q1) to[out=0,in=120] node[midway,above] {$\bra{m}$} (q4);
			\draw[->] (q2) to[out=0,in=-120] node[midway,below] {$\bra{m}$} (q4);
		\end{tikzpicture}
		\longmapsto
		\begin{tikzpicture}[baseline=-3pt,scale=0.9,font=\small]
		\node[myQV] (q1) at (-0.8,0.8) {};
		\node[myQV] (q2) at (-0.8,-0.8) {};
		\node[myQV] (q3) at (0,0) {};
		\node[myMQV] (q4) at (1.3,0) {};
		\draw[<-] (q3) -- (q1);
		\draw[->] (q3) -- (q2);
		\draw[<-] (q1) to[out=0,in=120] node[midway,above] {$\bra{m}$} (q4);
		\draw[->] (q2) to[out=0,in=-120] node[midway,below] {$\bra{m}$} (q4);
		\end{tikzpicture}
		\longmapsto
		\begin{tikzpicture}[baseline=-3pt,scale=0.9,font=\small]
		\node[myQV] (q1) at (-0.8,0.8) {};
		\node[myQV] (q2) at (-0.8,-0.8) {};
		\node[myQV] (q3) at (0,0) {};
		\node[myQV] (q4) at (1.3,0) {};
		\draw[<-] (q3) -- (q1);
		\draw[->] (q3) -- (q2);
		\draw[->] (q2) -- node[midway,left=-1pt] {$\bra{m}^2$} (q1);
		\draw[->] (q1) to[out=0,in=120] node[midway,above] {$\bra{m}$} (q4);
		\draw[<-] (q2) to[out=0,in=-120] node[midway,below] {$\bra{m}$} (q4);
		\end{tikzpicture}
		\]
		
		For $\mathcal H_4$, the following mutations give the weight $\bra{5}^2>2$.
		
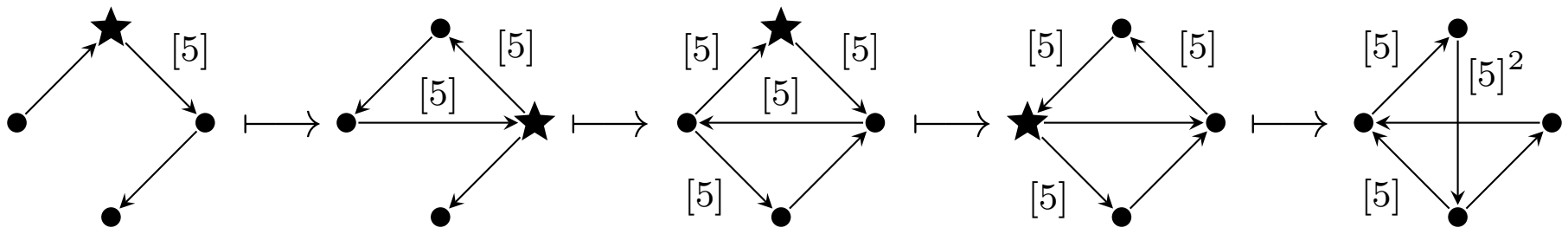
\[
		\begin{tikzpicture}[baseline=-3pt,scale=0.72,font=\scriptsize]
			\node[myQV] (q1) at (-1,0) {};
			\node[myMQV] (q2) at (0,1) {};
			\node[myQV] (q3) at (1,0) {};
			\node[myQV] (q4) at (0,-1) {};
			\draw[->] (q1) -- (q2);
			\draw[->] (q2) -- node[midway,above right=-2pt] {$\bra{5}$} (q3);
			\draw[->] (q3) -- (q4);
		\end{tikzpicture}
		\longmapsto
		\begin{tikzpicture}[baseline=-3pt,scale=0.72,font=\scriptsize]
			\node[myQV] (q1) at (-1,0) {};
			\node[myQV] (q2) at (0,1) {};
			\node[myMQV] (q3) at (1,0) {};
			\node[myQV] (q4) at (0,-1) {};
			\draw[->] (q1) -- node[midway,above=-2pt] {$\bra{5}$} (q3);
			\draw[->] (q2) -- (q1);
			\draw[->] (q3) -- node[midway,above right=-2pt] {$\bra{5}$} (q2);
			\draw[->] (q3) -- (q4);
		\end{tikzpicture}
		\longmapsto
		\begin{tikzpicture}[baseline=-3pt,scale=0.72,font=\scriptsize]
			\node[myQV] (q1) at (-1,0) {};
			\node[myMQV] (q2) at (0,1) {};
			\node[myQV] (q3) at (1,0) {};
			\node[myQV] (q4) at (0,-1) {};
			\draw[->] (q1) -- node[midway,above left=-2pt] {$\bra{5}$} (q2);
			\draw[->] (q1) -- node[midway,below left=-2pt] {$\bra{5}$} (q4);
			\draw[->] (q2) -- node[midway,above right=-2pt] {$\bra{5}$} (q3);
			\draw[->] (q3) -- node[midway,above=-2pt] {$\bra{5}$} (q1);
			\draw[->] (q4) -- (q3);
		\end{tikzpicture}
		\longmapsto
		\begin{tikzpicture}[baseline=-3pt,scale=0.72,font=\scriptsize]
			\node[myMQV] (q1) at (-1,0) {};
			\node[myQV] (q2) at (0,1) {};
			\node[myQV] (q3) at (1,0) {};
			\node[myQV] (q4) at (0,-1) {};
			\draw[->] (q1) -- (q3);
			\draw[->] (q1) -- node[midway,below left=-2pt] {$\bra{5}$} (q4);
			\draw[->] (q2) -- node[midway,above left=-2pt] {$\bra{5}$} (q1);
			\draw[->] (q3) -- node[midway,above right=-2pt] {$\bra{5}$} (q2);
			\draw[->] (q4) -- (q3);
		\end{tikzpicture}
		\longmapsto
		\begin{tikzpicture}[baseline=-3pt,scale=0.72,font=\scriptsize]
			\node[myQV] (q1) at (-1,0) {};
			\node[myQV] (q2) at (0,1) {};
			\node[myQV] (q3) at (1,0) {};
			\node[myQV] (q4) at (0,-1) {};
			\draw[->] (q1) -- node[midway,above left=-2pt] {$\bra{5}$} (q2);
			\draw[->] (q2) -- node[midway,right=-2pt,pos=.2] {$\bra{5}^2$} (q4);
			\draw[->] (q3) -- (q1);
			\draw[->] (q4) -- node[midway,below left=-2pt] {$\bra{5}$} (q1);
			\draw[->] (q4) -- (q3);
		\end{tikzpicture}
		\]
		
		For $Z_5(3,5)$, the following mutations give a quiver whose bold arrows form a full subquiver that is a reorientation of $\mathcal H_4$. Thus the above applies.
		
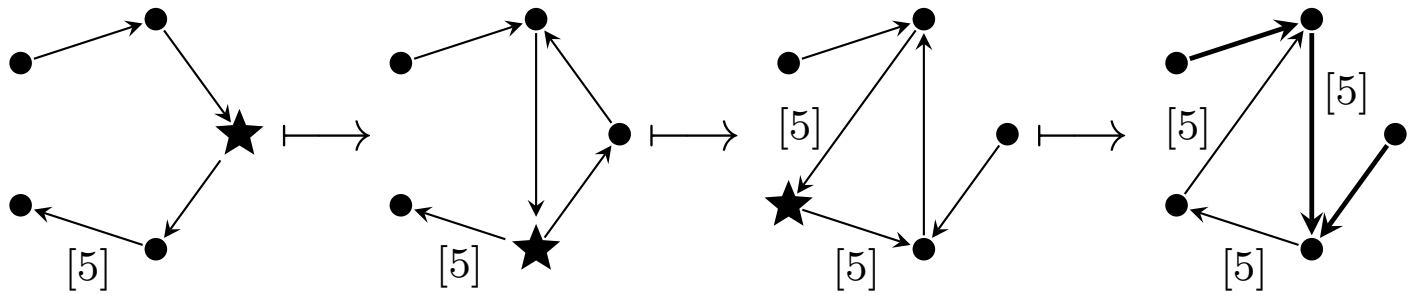
\[
		\begin{tikzpicture}[baseline=-3pt,scale=0.8,font=\scriptsize]
			\node[myQV] (q1) at (-0.80902,0.58779) {};
			\node[myQV] (q2) at (0.30902,0.95106) {};
			\node[myMQV] (q3) at (1.0,0.0) {};
			\node[myQV] (q4) at (0.30902,-0.95106) {};
			\node[myQV] (q5) at (-0.80902,-0.58779) {};
			\draw[->] (q1) -- (q2);
			\draw[->] (q2) -- (q3);
			\draw[->] (q3) -- (q4);
			\draw[->] (q4) -- node[midway,below] {$\bra{5}$} (q5);
		\end{tikzpicture}
		\longmapsto
		\begin{tikzpicture}[baseline=-3pt,scale=0.8,font=\scriptsize]
			\node[myQV] (q1) at (-0.80902,0.58779) {};
			\node[myQV] (q2) at (0.30902,0.95106) {};
			\node[myQV] (q3) at (1.0,0.0) {};
			\node[myMQV] (q4) at (0.30902,-0.95106) {};
			\node[myQV] (q5) at (-0.80902,-0.58779) {};
			\draw[->] (q1) -- (q2);
			\draw[->] (q2) -- (q4);
			\draw[->] (q3) -- (q2);
			\draw[->] (q4) -- (q3);
			\draw[->] (q4) -- node[midway,below] {$\bra{5}$} (q5);
		\end{tikzpicture}
		\longmapsto
		\begin{tikzpicture}[baseline=-3pt,scale=0.8,font=\scriptsize]
			\node[myQV] (q1) at (-0.80902,0.58779) {};
			\node[myQV] (q2) at (0.30902,0.95106) {};
			\node[myQV] (q3) at (1.0,0.0) {};
			\node[myQV] (q4) at (0.30902,-0.95106) {};
			\node[myMQV] (q5) at (-0.80902,-0.58779) {};
			\draw[->] (q1) -- (q2);
			\draw[->] (q2) -- node[midway,left,pos=.6] {$\bra{5}$} (q5);
			\draw[->] (q3) -- (q4);
			\draw[->] (q4) -- (q2);
			\draw[->] (q5) -- node[midway,below] {$\bra{5}$} (q4);
		\end{tikzpicture}
		\longmapsto
		\begin{tikzpicture}[baseline=-3pt,scale=0.8,font=\scriptsize]
			\node[myQV] (q1) at (-0.80902,0.58779) {};
			\node[myQV] (q2) at (0.30902,0.95106) {};
			\node[myQV] (q3) at (1.0,0.0) {};
			\node[myQV] (q4) at (0.30902,-0.95106) {};
			\node[myQV] (q5) at (-0.80902,-0.58779) {};
			\draw[->,line width=.9pt] (q1) -- (q2);
			\draw[->,line width=.9pt] (q2) -- node[midway,right=-2pt,pos=.3] {$\bra{5}$} (q4);
			\draw[->,line width=.9pt] (q3) -- (q4);
			\draw[->] (q4) -- node[midway,below] {$\bra{5}$} (q5);
			\draw[->] (q5) -- node[midway,left,pos=.4] {$\bra{5}$} (q2);
		\end{tikzpicture}
		\]
	\end{proof}
	
	\begin{lemma}\label{lem: 2 infiniteness of Y and Z}
		The following statements hold.
		\begin{enumerate}[label=\textup{(\arabic*)}]
			\item $Y_n(m)$ is $2$-infinite for $m,n\in\mathbb Z_{\geq4}$.
			\item Let $m,m',n\in\mathbb Z_{\geq3}$. Then $Z_n(m,m')$ is $2$-infinite unless
			\[
			\{m,m'\}=\{3\},\{3,4\},
			\quad\text{or}\quad
			\{m,m'\}=\{3,5\}\ \text{with }n\in\{3,4\}.
			\]
		\end{enumerate}
	\end{lemma}
	
	\begin{proof}
		We use the following shortening operation. At a vertex $k$ of degree $2$ in a tree, Lemma~\ref{lem: chain reorientations} allows us to orient the two incident arrows as shown below. Mutating at $k$ and then deleting it gives
		\[
		i\xrightarrow{a}k\xrightarrow{b}j
		\quad\xrightarrow{\ \mu_k,\ \text{delete }k\ }\quad
		i\xrightarrow{ab}j.
		\]
		In particular, if one of $a,b$ is $1$, then the other weight is preserved.
		
		For (1), this operation shortens the long branch of $Y_n(m)$ to give $Y_4(m)$, which is $2$-infinite by Lemma~\ref{lem: special 2 infinite quivers}.
		
		For (2), we may assume that $m\leq m'$. First suppose that $m,m'\geq4$. Shortening gives a three-vertex path with weights $\bra{m}$ and $\bra{m'}$. Mutation at its middle vertex gives the weight $\bra{m}\bra{m'}\geq2$. Next suppose that $m=3$ and $m'\geq6$. Shortening gives a path with weights $1$ and $\bra{m'}$. Two mutations give the weight $\bra{m'}^2-1\geq2$, as in the proof of Lemma~\ref{lem: rank 3 tree case}. Finally, when $(m,m')=(3,5)$ and $n\geq5$, shortening gives $Z_5(3,5)$, which is $2$-infinite by Lemma~\ref{lem: special 2 infinite quivers}. These are all cases not listed in the statement.
	\end{proof}
	We use the following known fact.
	\begin{lemma}[{\cite[proof of Proposition~9.3]{FZ03a}}]\label{lem:affine-obstructions}
		Every orientation of the following weighted tree is $2$-infinite:
		\[
		\widetilde F_4:\quad
		\begin{tikzpicture}[baseline=-3pt]
			\foreach \x in {0,...,4}{
				\node[myQV] (v\x) at (\x,0) {};
			}
			\draw (v0) -- (v1) -- (v2);
			\draw (v2) -- node[midway,above] {$\bra{4}$} (v3);
			\draw (v3) -- (v4);
		\end{tikzpicture}
		\]
	\end{lemma}
	
	Here we use the square root of the diagram weight in \cite{FZ03a}.
	
	\begin{proof}[Proof of Proposition~\ref{prop: tree Coxeter quivers and 2-finite}]
		The ``if'' implication is \eqref{eq:Cox-2fin-DB2}. Suppose that $Q \in \D$ is $2$-finite. Its weights are $\bra m$ with $m\geq3$. The assertion is immediate in ranks $1$ and $2$, so assume that $Q$ has rank at least $3$.
		
		If all weights are $1$, then the claim follows from \cite[Proposition~9.3]{FZ03a}. We may therefore assume that $Q$ has a non-unit weight. If $Q$ had more than one non-unit edge, then it would contain a full subquiver $Z_k(m,m')$ with $m,m'\geq4$ and $k \geq 3$, contrary to Lemma~\ref{lem: 2 infiniteness of Y and Z}(2). Hence $Q$ has exactly one non-unit edge, of weight $\bra{m}$ with $m\geq4$. If $Q$ had a branching vertex, then it would contain a full subquiver $Y_k(m)$, contrary to Lemma~\ref{lem: 2 infiniteness of Y and Z}(1). Thus $Q$ is a path.
		
		If $m=4$, then a terminal non-unit edge gives type $B_n$. If the edge is internal, then a path with at least five vertices would contain $\widetilde F_4$, contrary to Lemma~\ref{lem:affine-obstructions}. Thus an internal non-unit edge gives type $F_4$. If $m=5$, then Lemma~\ref{lem: special 2 infinite quivers} excludes $\mathcal H_4$, so the non-unit edge must be terminal. Lemma~\ref{lem: special 2 infinite quivers} excludes $Z_5(3,5)$ and implies that $Q$ has rank at most $4$. This gives $H_3$ or $H_4$. If $m\geq6$, then $Q$ contains $Z_3(3,m)$, contrary to Lemma~\ref{lem: 2 infiniteness of Y and Z}(2).
	\end{proof}

\subsection{Classification of chordless cycles}
	We next treat chordless cycles and prove the following classification.
	
	\begin{proposition}[cf.~{\cite[Proposition~9.7]{FZ03a}}]\label{prop: cyclically oriented lemma}
		Let $Q\in\D$ be a quiver whose underlying weighted graph is an $n$-cycle. Then $Q$ is $2$-finite if and only if it is a directed cycle and satisfies one of the following conditions:
    \begin{enumerate}[label=\textup{(\arabic*)}]
      \item all weights are $1$, giving mutation type $D_n$ \textup{(with $D_3=A_3$)};
      \item $n=3$ and the weights are $(1,\bra{4},\bra{4})$, $(1,\bra{5},\bra{5})$, or $(\bra{5},\bra{5},\bra{5})$, giving types $B_3$, $H_3$, and $H_3$, respectively, as listed in Table~\ref{tab: 2-finite of rank 3};
      \item $n=4$ and the consecutive weights are $({\bra{4}},1,{\bra{4}},1)$ or $(\bra{5},1,\bra{5},1)$, giving types $F_4$ and $H_4$, respectively (see \eqref{eq:typeF_4cycle} and \eqref{eq:typeH_4cycle}).
    \end{enumerate}
	\end{proposition}
	
	To prove Proposition~\ref{prop: cyclically oriented lemma}, we prepare two lemmas.
	
	\begin{lemma}\label{lem: cycle restriction from D}
		Let $Q$ be a directed cycle with two distinguished arrows of weights $\bra m$ and $\bra{m'}$, where $m>m'\geq3$. If all other arrows have weight $1$, then $Q\notin\D$.
	\end{lemma}
	
	\begin{proof}
		For a triangle, mutation at the vertex between the arrows of weights $1$ and $\bra{m'}$ gives the weight $\bra m-\bra{m'}$. Since $0<\bra m-\bra{m'}<1$, $Q\notin\D$.
		
		For a longer cycle, choose an endpoint $k$ of an arrow of weight $1$. Then $\mu_k(Q)\setminus\{k\}$ is a shorter directed cycle with weights $\bra m$, $\bra{m'}$, and $1$ on all remaining arrows. The assertion follows by induction on the length of the cycle.
	\end{proof}
	
	\begin{lemma}[cf.~{\cite[Proposition~9.7]{FZ03a}}]\label{lem: cyclic orientation of 2 finite quivers}
		Every chordless cycle in a $2$-finite quiver $Q\in\D$ is a directed cycle.
	\end{lemma}
	
	\begin{proof}
		We may assume that the underlying weighted graph of $Q$ is an $n$-cycle. We use induction on $n$. The case $n=3$ follows from Proposition~\ref{prop: 2-finte rank 3 quivers}.
		
		Suppose that $n\geq4$. If a vertex $k$ has one incoming and one outgoing arrow, then $\mu_k(Q)\setminus\{k\}$ is a chordless $(n-1)$-cycle. By induction it is directed, and hence so is $Q$. Otherwise, every vertex is a source or sink. Mutation at any vertex gives a cycle in which its neighbors have one incoming and one outgoing arrow. This contradicts the preceding argument.
	\end{proof}
	
	Note that the assumption $Q\in\D$ in Lemma~\ref{lem: cyclic orientation of 2 finite quivers} cannot be omitted: mutating the $2$-finite quiver in Example~\ref{ex:2fin-non-Cox} at any vertex gives an acyclic triangle.
	
	\begin{proof}[Proof of Proposition~\ref{prop: cyclically oriented lemma}]
		We first check the quivers in the list. For (1), Lemma~\ref{lem: mutation equivalence between S and T} shows that $S_{1,1,1}^{n-3}$ is mutation-equivalent to $T_{1,1,n-3}$, which has type $D_n$. The triangles in (2) were treated in Proposition~\ref{prop: 2-finte rank 3 quivers}. For (3), the following mutations, applied in the indicated order, give finite Coxeter quivers of types $F_4$ and $H_4$.
		\begin{gather}
			\begin{tikzpicture}[baseline=-3pt,scale=0.75]
				\node[myQV,label={above left:$1$}] (q1) at (-0.8,0.65) {};
				\node[myQV,label={above right:$2$}] (q2) at (0.8,0.65) {};
				\node[myQV,label={below right:$3$}] (q3) at (0.8,-0.65) {};
				\node[myQV,label={below left:$4$}] (q4) at (-0.8,-0.65) {};
				\draw[->] (q1) -- node[midway,above] {$\bra{4}$} (q2);
				\draw[->] (q2) -- (q3);
				\draw[->] (q3) -- node[midway,below] {$\bra{4}$} (q4);
				\draw[->] (q4) -- (q1);
			\end{tikzpicture}
			\quad\xrightarrow{\ \mu_{1},\mu_{2},\mu_{1}\ }\quad
			\begin{tikzpicture}[baseline=-3pt,scale=0.9]
				\node[myQV,label={below:$4$}] (q3) at (0.0,0) {};
				\node[myQV,label={below:$1$}] (q2) at (1.25,0) {};
				\node[myQV,label={below:$2$}] (q1) at (2.5,0) {};
				\node[myQV,label={below:$3$}] (q4) at (3.75,0) {};
				\draw[<-] (q1) -- (q4);
				\draw[<-] (q2) -- node[midway,above] {$\bra{4}$} (q1);
				\draw[<-] (q3) -- (q2);
			\end{tikzpicture}
      \label{eq:typeF_4cycle}
      \\[12pt]
      \begin{tikzpicture}[baseline=-3pt,scale=0.75]
				\node[myQV,label={above left:$1$}] (q1) at (-0.8,0.65) {};
				\node[myQV,label={above right:$2$}] (q2) at (0.8,0.65) {};
				\node[myQV,label={below right:$3$}] (q3) at (0.8,-0.65) {};
				\node[myQV,label={below left:$4$}] (q4) at (-0.8,-0.65) {};
				\draw[->] (q1) -- node[midway,above] {$\bra{5}$} (q2);
				\draw[->] (q2) -- (q3);
				\draw[->] (q3) -- node[midway,below] {$\bra{5}$} (q4);
				\draw[->] (q4) -- (q1);
			\end{tikzpicture}
			\quad\xrightarrow{\ \mu_{1},\mu_{2},\mu_{4},\mu_{1}\ }\quad
			\begin{tikzpicture}[baseline=-3pt,scale=0.9]
				\node[myQV,label={below:$4$}] (q4) at (0.0,0) {};
				\node[myQV,label={below:$1$}] (q1) at (1.25,0) {};
				\node[myQV,label={below:$2$}] (q2) at (2.5,0) {};
				\node[myQV,label={below:$3$}] (q3) at (3.75,0) {};
				\draw[->] (q1) -- (q2);
				\draw[->] (q3) -- (q2);
				\draw[->] (q4) -- node[midway,above] {$\bra{5}$} (q1);
			\end{tikzpicture}
      \label{eq:typeH_4cycle}
		\end{gather}
		Thus \eqref{eq:Cox-2fin-DB2} shows that every quiver in the list is $2$-finite.
		
		Conversely, suppose that $Q$ is $2$-finite. It is a directed cycle by Lemma~\ref{lem: cyclic orientation of 2 finite quivers}. The case $n=3$ follows from Proposition~\ref{prop: 2-finte rank 3 quivers}, so assume that $n\geq4$. If all weights are $1$, then $Q$ is in case~(1). Otherwise, Lemma~\ref{lem: cycle restriction from D} implies that $Q$ has at least two non-unit edges. Two non-unit edges cannot be consecutive, since they would give a full subquiver $Z_3(m,m')$ with $m,m'\geq4$, contrary to Lemma~\ref{lem: 2 infiniteness of Y and Z}(2).
		
		If $n\geq5$, then $Q$ would contain a full subquiver $Z_k(m,m')$ with $m,m'\geq4$, contrary to Lemma~\ref{lem: 2 infiniteness of Y and Z}(2). Thus $n=4$. Exactly two opposite edges have non-unit weights and their weights must be equal by Lemma~\ref{lem: cycle restriction from D}. By Lemma~\ref{lem: rank 3 tree case}, applied to three consecutive vertices, the weight is $\bra{4}$ or $\bra{5}$.
	\end{proof}

\subsection{Preparations for the induction}
	We use the following lemmas for the induction in the proof of Theorem~\ref{thm:2fin-Cox}. In the diagrams below, dotted arrows denote directed paths, and weights and edge orientations are omitted when they are not needed.
	
	\begin{lemma}\label{lem:cycle-attachment}
		Let $Q\in\D$ be $2$-finite, $\mathcal C\subseteq Q$ a chordless cycle with at least four vertices, and $u$ a vertex outside $\mathcal C$. Then $u$ has at most two neighbors in $\mathcal C$. If it has two, then they are adjacent in $\mathcal C$.
	\end{lemma}
	
	\begin{proof}
		By Lemma~\ref{lem: cyclic orientation of 2 finite quivers}, every chordless cycle in $Q$ is directed. If $u$ has exactly two neighbors in $\mathcal C$ and they are nonadjacent, then the two chordless cycles through $u$ require opposite orientations of each edge at $u$, as in the diagram on the left below. If $u$ has at least three neighbors, then we consider three successive ones along $\mathcal C$. The two chordless cycles as in the diagram on the right below require opposite orientations of their common edge at $u$.
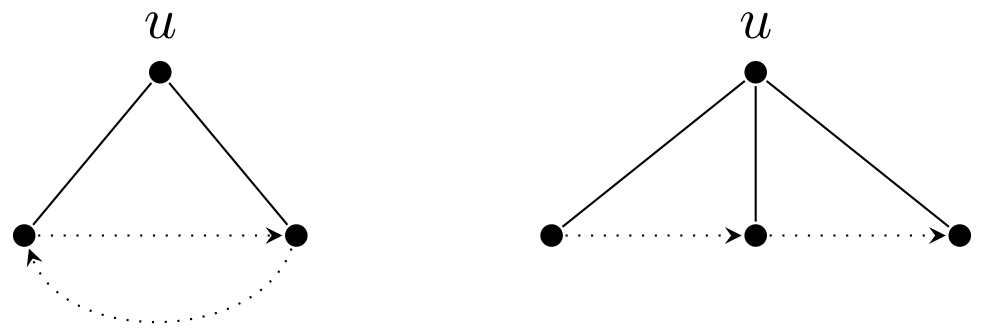
		\[
		\begin{tikzpicture}[baseline=-3pt,scale=0.9]
			\node[myQV] (x) at (-1,0) {};
			\node[myQV] (y) at (1,0) {};
			\node[myQV,label=above:$u$] (u) at (0,1.2) {};
			\draw[->,dotted] (y) to[bend left=70] (x);
			\draw[->,dotted] (x) to (y);
			\draw (x)--(u)--(y);
		\end{tikzpicture}
		\qquad\qquad
		\begin{tikzpicture}[baseline=-3pt,scale=0.9]
			\node[myQV] (x) at (-1.5,0) {};
			\node[myQV] (y) at (1.5,0) {};
			\node[myQV] (z) at (0,0) {};
			\node[myQV,label=above:$u$] (u) at (0,1.2) {};
			\draw[->,dotted] (x) to (z);
			\draw[->,dotted] (z) to (y);
			\draw (x)--(u)--(y);
			\draw (u)--(z);
		\end{tikzpicture}
		\]
		Both cases give a contradiction.
	\end{proof}
	
	\begin{lemma}\label{lem: eliminating 3-cycle lemma}
		Let $\mathcal C$ be a directed triangle with weights $(1,a,a)$, where $a>0$. Attach two disjoint paths to distinct vertices of $\mathcal C$, each meeting $\mathcal C$ only at an endpoint. Either path may have length zero. If all arrows of at least one path have weight $1$, then the resulting quiver is mutation-equivalent to a tree.
	\end{lemma}
	
	\begin{proof}
		If an endpoint of the arrow of weight $1$ in $\mathcal C$ has no path of positive length attached to it, then the following mutation gives a tree.
		
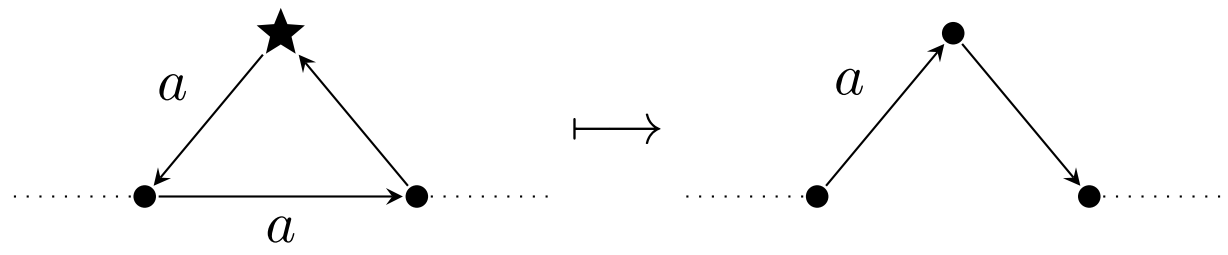
\[
		\begin{tikzpicture}[baseline=10pt,scale=0.9]
			\node[myQV] (qa) at (-1,0) {};
			\node[myQV] (qb) at (1,0) {};
			\node[myMQV] (qu) at (0,1.2) {};
			\draw[->] (qa) -- node[midway,below] {$a$} (qb);
			\draw[->] (qb) -- (qu);
			\draw[->] (qu) -- node[midway,above left] {$a$} (qa);
			\draw[dotted] (qa)--(-2,0);
			\draw[dotted] (qb)--(2,0);
		\end{tikzpicture}
		\longmapsto
		\begin{tikzpicture}[baseline=10pt,scale=0.9]
			\node[myQV] (qa) at (-1,0) {};
			\node[myQV] (qb) at (1,0) {};
			\node[myQV] (qu) at (0,1.2) {};
			\draw[<-] (qb) -- (qu);
			\draw[<-] (qu) -- node[midway,above left] {$a$} (qa);
			\draw[dotted] (qa)--(-2,0);
			\draw[dotted] (qb)--(2,0);
		\end{tikzpicture}
		\]
		Otherwise, we consider an attached path with only unit weights. By Lemma~\ref{lem: chain reorientations} and the obvious symmetry, we may orient the resulting quiver as on the left below. The following mutation preserves the form of the quiver and shortens the attached path by one edge:
		
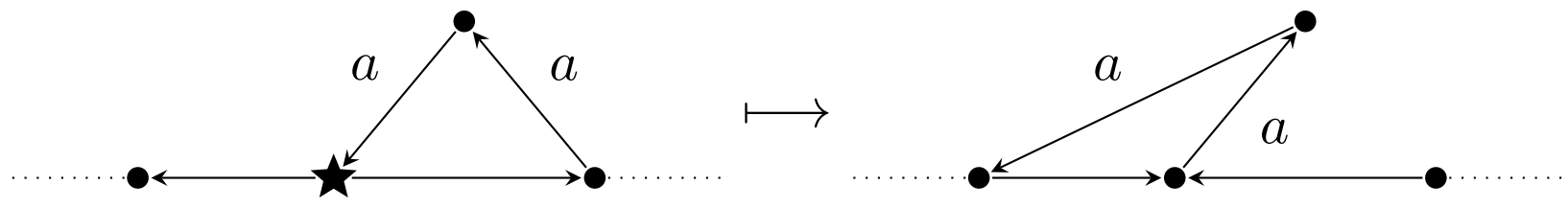
\[
		\begin{tikzpicture}[baseline=10pt,scale=0.9]
			\node[myMQV] (qa) at (-1,0) {};
			\node[myQV] (qb) at (1,0) {};
			\node[myQV] (qc) at (-2.5,0) {};
			\node[myQV] (qu) at (0,1.2) {};
			\draw[->] (qa) -- (qb);
			\draw[->] (qb) -- node[midway,above right] {$a$} (qu);
			\draw[->] (qu) -- node[midway,above left] {$a$} (qa);
			\draw[->] (qa) -- (qc);
			\draw[dotted] (qc)--(-3.5,0);
			\draw[dotted] (qb)--(2,0);
		\end{tikzpicture}
		\longmapsto
		\begin{tikzpicture}[baseline=10pt,scale=0.9]
			\node[myQV] (qa) at (-1,0) {};
			\node[myQV] (qb) at (1,0) {};
			\node[myQV] (qc) at (-2.5,0) {};
			\node[myQV] (qu) at (0,1.2) {};
			\draw[<-] (qa) -- (qb);
			\draw[->] (qu) -- node[midway,above left] {$a$} (qc);
			\draw[<-] (qu) -- node[midway,below right] {$a$} (qa);
			\draw[<-] (qa) -- (qc);
			\draw[dotted] (qc)--(-3.5,0);
			\draw[dotted] (qb)--(2,0);
		\end{tikzpicture}
		\]
		The assertion follows by induction on the length of this path.
	\end{proof}

	\begin{lemma}\label{lem: maximality of F4 and H4}
		Let $Q\in\D$ be connected. If it contains a directed $4$-cycle of mutation type $F_4$ or $H_4$ as a proper full subquiver, then $Q$ is $2$-infinite.
	\end{lemma}
	
	\begin{proof}
		Suppose that $Q$ is $2$-finite. Let $\mathcal C$ be the given cycle with weights $(\bra{m},1,\bra{m},1)$ in cyclic order, where $m\in\{4,5\}$. By connectedness, we may restrict to the full subquiver on $\mathcal C$ and a vertex $u$ adjacent to it. By Lemma~\ref{lem:cycle-attachment}, $u$ has one neighbor or two adjacent neighbors in $\mathcal C$.
		
		If $u$ has one neighbor, then $Q$ contains a full subquiver $Y_4(m)$ or $Z_3(m,m')$ with $m'\geq4$, contrary to Lemma~\ref{lem: 2 infiniteness of Y and Z}. If $u$ has two neighbors, then Proposition~\ref{prop: 2-finte rank 3 quivers} shows that mutation at $u$ either reverses or deletes the edge joining them. In the first case, $\mu_u(Q)\setminus\{u\}$ is a cycle that is not directed, contrary to Lemma~\ref{lem: cyclic orientation of 2 finite quivers}. In the second case, $\mu_u(Q)$ is a chordless $5$-cycle with a non-unit edge, contrary to Proposition~\ref{prop: cyclically oriented lemma}.
	\end{proof}

\subsection{Completion of the proof}
	When all weights are $1$, the exchange matrix is integral and skew-symmetric. In this case, $2$-finiteness in our sense is equivalent to $2$-finiteness in the sense of \cite{FZ03a}. Thus \cite[Theorem~8.6]{FZ03a} gives the following statement.
	\begin{equation}\label{eq:unit-weight-classification}
		\begin{gathered}
			\text{a connected $2$-finite quiver with all weights $1$ is of finite Coxeter type.}
		\end{gathered}
	\end{equation}
	By using this fact, we complete the proof of Theorem~\ref{thm:2fin-Cox}.
	\begin{proof}[Proof of Theorem~\ref{thm:2fin-Cox}]
		The ``if'' implication is \eqref{eq:Cox-2fin-DB2}. Let $Q\in\D$ be connected and $2$-finite. We prove that $Q$ is of finite Coxeter type by induction on its rank. The cases of ranks $1$ and $2$ are immediate, and rank $3$ follows from Proposition~\ref{prop: 2-finte rank 3 quivers}.
		
		Suppose that $Q$ has rank $n+1$, where $n\geq3$. Choose a vertex $u$ such that $Q'=Q\setminus\{u\}$ is connected. By induction, $Q'$ is of finite Coxeter type. Since mutations of $Q'$ extend to $Q$ by \eqref{eq:restriction-mutation}, we may replace $Q'$ by any quiver in its mutation class. Whenever we obtain a tree, Proposition~\ref{prop: tree Coxeter quivers and 2-finite} completes the argument.
		
		\medskip
		\noindent\emph{Case 1: $Q'$ has type $A_n$, $B_n$, $F_4$, $H_3$, or $H_4$.}
		Choose $Q'$ to be the corresponding directed path. By the same argument as in the proof of Lemma~\ref{lem:cycle-attachment}, $u$ has at most two neighbors in $Q'$. If it has one neighbor, then $Q$ is a tree. Otherwise, $Q$ has a unique chordless cycle $\mathcal C$, with paths attached at two of its vertices. We consider the possibilities for $\mathcal C$ given by Proposition~\ref{prop: cyclically oriented lemma}.
			
		Suppose first that $\mathcal C$ is a triangle. If its weights are $(1,\bra{m},\bra{m})$ with $m\in\{3,4,5\}$, then Lemma~\ref{lem: eliminating 3-cycle lemma} implies that $Q$ is mutation-equivalent to a tree since $Q'$ has at most one non-unit edge. If all its weights are $\bra{5}$, then $\mu_u(Q)$ is the above case with weights $(1,\bra{5},\bra{5})$.
		
		If $\mathcal C$ has at least four vertices and all its weights are $1$, then any non-unit edge outside $\mathcal C$ would give a full subquiver $Y_k(m)$ with $m\geq4$, contrary to Lemma~\ref{lem: 2 infiniteness of Y and Z}(1). Hence all weights of $Q$ are $1$, and \eqref{eq:unit-weight-classification} applies.
		
		Finally, if $\mathcal C$ is a four-cycle of mutation type $F_4$ or $H_4$, then Lemma~\ref{lem: maximality of F4 and H4} forces $Q=\mathcal C$, whose mutation type is already known.
		
		\medskip
		\noindent\emph{Case 2: $Q'$ has type $D_n$, with $n\geq4$.}
		By Proposition~\ref{prop: cyclically oriented lemma}(1), choose $Q'$ as the directed $n$-cycle with all weights $1$. By Lemma~\ref{lem:cycle-attachment}, $u$ has one neighbor or two adjacent neighbors in $Q'$. In the first case, $Q$ is a tree reorientation of $S_{2,1,1}^{n-3}(m)$ for some $m\geq3$, so Lemma~\ref{lem: mutation equivalence between S and T} gives a tree. In the second case, Proposition~\ref{prop: 2-finte rank 3 quivers} shows that the triangle containing $u$ is directed, and its weights are $(1,\bra{m},\bra{m})$ with $m\in\{3,4,5\}$. Mutation at $u$ gives
		
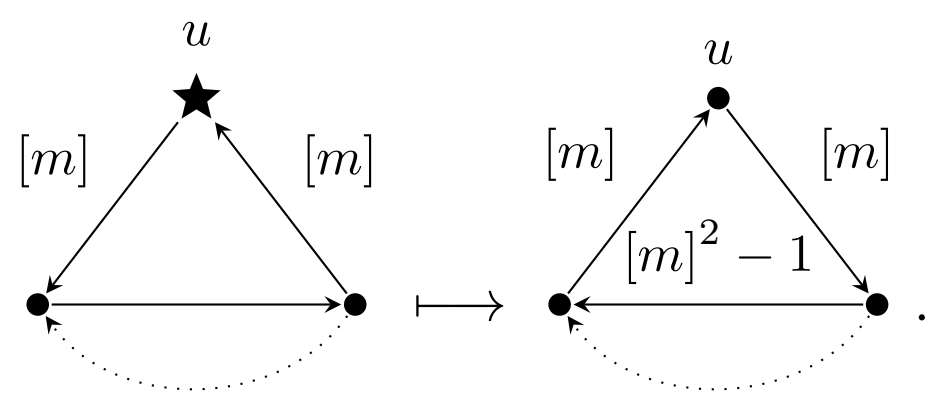
\[
		\begin{tikzpicture}[baseline=-3pt,scale=1.05,font=\small]
			\node[myQV] (qx) at (-1,0) {};
			\node[myQV] (qy) at (1,0) {};
			\node[myMQV,label={[font=\small]above:$u$}] (qu) at (0,1.3) {};
			\draw[->] (qx) -- (qy);
			\draw[->] (qy) -- node[midway,above right] {$\bra{m}$} (qu);
			\draw[->] (qu) -- node[midway,above left] {$\bra{m}$} (qx);
			\draw[->,dotted] (qy) to[bend left=55] (qx);
		\end{tikzpicture}
		\longmapsto
		\begin{tikzpicture}[baseline=-3pt,scale=1.05,font=\small]
			\node[myQV] (qx) at (-1,0) {};
			\node[myQV] (qy) at (1,0) {};
			\node[myQV,label={[font=\small]above:$u$}] (qu) at (0,1.3) {};
			\draw[->] (qx) -- node[midway,above left] {$\bra{m}$} (qu);
			\draw[->] (qu) -- node[midway,above right] {$\bra{m}$} (qy);
			\draw[->] (qy) -- node[midway,above] {$\bra{m}^2-1$} (qx);
			\draw[->,dotted] (qy) to[bend left=55] (qx);
		\end{tikzpicture}.
		\]
		If $m>3$, then $\mu_u(Q)\setminus\{u\}$ is a cycle that is not directed, contrary to Lemma~\ref{lem: cyclic orientation of 2 finite quivers}. Hence $m=3$, and $\mu_u(Q)$ is the directed $(n+1)$-cycle with all weights $1$, of type $D_{n+1}$ by Proposition~\ref{prop: cyclically oriented lemma}(1).
		
		\medskip
		\noindent\emph{Case 3: $Q'$ has type $E_n$, where $n=6,7,8$.}
		By Lemma~\ref{lem: mutation equivalence between S and T}, take $Q'=S_{2,1,1}^{n-4}$, consisting of a directed $(n-1)$-cycle $\mathcal C$ and a leaf $v_1$ attached at $v_0$, with all weights $1$:
		
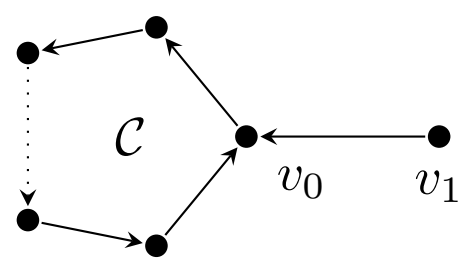
\[
		\begin{tikzpicture}[baseline=(current bounding box.center),scale=.85,>=stealth,font=\small]
			\node[myQV,label=below right:$v_0$] (v0) at (0,0) {};
			\node[myQV] (c1) at (-.7,.85) {};
			\node[myQV] (c2) at (-1.7,.65) {};
			\node[myQV] (c3) at (-1.7,-.65) {};
			\node[myQV] (c4) at (-.7,-.85) {};
			\draw[->] (v0)--(c1);
			\draw[->] (c1)--(c2);
			\draw[->,dotted] (c2)--(c3);
			\draw[->] (c3)--(c4);
			\draw[->] (c4)--(v0);
			\node at (-.9,0) {$\mathcal C$};
			\node[myQV,label=below:$v_1$] (v1) at (1.5,0) {};
			\draw[->] (v1)--(v0);
		\end{tikzpicture}
		\]
		By Lemma~\ref{lem:cycle-attachment}, $u$ has at most two neighbors in $\mathcal C$, and two such neighbors must be adjacent.
		
		If $u$ has no neighbor in $\mathcal C$, then it is adjacent only to $v_1$. Thus $Q$ is a tree reorientation of $S_{3,1,1}^{n-4}(3,m)$ for some $m\geq3$, and Lemma~\ref{lem: mutation equivalence between S and T} gives a tree.
		
		Suppose that $u$ has exactly one neighbor $x$ in $\mathcal C$. If $u$ is not adjacent to $v_1$, then the edge between $u$ and $x$ must have weight $1$: otherwise, $Q$ would contain a full subquiver $Y_4(m)$ with $m\geq4$, contrary to Lemma~\ref{lem: 2 infiniteness of Y and Z}(1). Thus all weights are $1$, and \eqref{eq:unit-weight-classification} applies. Now suppose that $u$ is also adjacent to $v_1$. If $x=v_0$, then Proposition~\ref{prop: 2-finte rank 3 quivers} gives the following mutation, with $m\in\{3,4,5\}$:
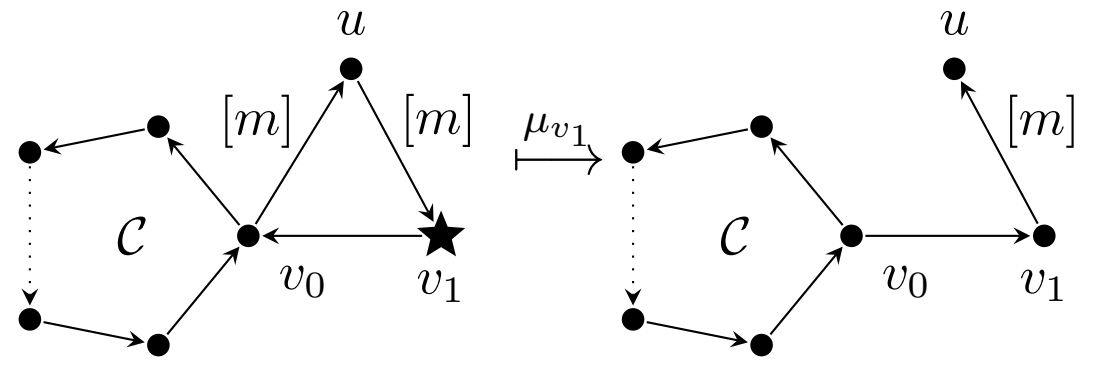
		\[
		\begin{tikzpicture}[baseline=(current bounding box.center),scale=.85,>=stealth,font=\small]
			\node[myQV,label=below right:$v_0$] (v0) at (0,0) {};
			\node[myQV] (c1) at (-.7,.85) {};
			\node[myQV] (c2) at (-1.7,.65) {};
			\node[myQV] (c3) at (-1.7,-.65) {};
			\node[myQV] (c4) at (-.7,-.85) {};
			\draw[->] (v0)--(c1);
			\draw[->] (c1)--(c2);
			\draw[->,dotted] (c2)--(c3);
			\draw[->] (c3)--(c4);
			\draw[->] (c4)--(v0);
			\node at (-.9,0) {$\mathcal C$};
			\node[myMQV,label=below:$v_1$] (v1) at (1.5,0) {};
			\node[myQV,label=above:$u$] (u) at (.8,1.3) {};
			\draw[->] (v1)--(v0);
			\draw[->] (v0)--node[above left=-3pt] {$\bra{m}$} (u);
			\draw[->] (u)--node[above right=-3pt] {$\bra{m}$} (v1);
		\end{tikzpicture}
		\overset{\mu_{v_1}}{\longmapsto}
		\begin{tikzpicture}[baseline=(current bounding box.center),scale=.85,>=stealth,font=\small]
			\node[myQV,label=below right:$v_0$] (v0) at (0,0) {};
			\node[myQV] (c1) at (-.7,.85) {};
			\node[myQV] (c2) at (-1.7,.65) {};
			\node[myQV] (c3) at (-1.7,-.65) {};
			\node[myQV] (c4) at (-.7,-.85) {};
			\draw[->] (v0)--(c1);
			\draw[->] (c1)--(c2);
			\draw[->,dotted] (c2)--(c3);
			\draw[->] (c3)--(c4);
			\draw[->] (c4)--(v0);
			\node at (-.9,0) {$\mathcal C$};
			\node[myQV,label=below:$v_1$] (v1) at (1.5,0) {};
			\node[myQV,label=above:$u$] (u) at (.8,1.3) {};
			\draw[->] (v0)--(v1);
			\draw[->] (v1)--node[above right=-3pt] {$\bra{m}$} (u);
		\end{tikzpicture}
		\]
		The result is a tree reorientation of $S_{3,1,1}^{n-4}(3,m)$, so Lemma~\ref{lem: mutation equivalence between S and T} gives a tree. If $x\neq v_0$, then the edges $v_1-u$ and $u-x$ belong to a chordless cycle of length at least $4$ whose other edges have weight $1$. Since these two edges are consecutive, their weights are $1$ by Proposition~\ref{prop: cyclically oriented lemma}. Thus \eqref{eq:unit-weight-classification} applies.
		
		Finally, suppose that $u$ has two neighbors in $\mathcal C$. The argument in Case~2 shows that both edges between $u$ and $\mathcal C$ have weight $1$. If $u$ is also adjacent to $v_1$, then this edge belongs to a chordless cycle whose other edges have weight $1$. By Proposition~\ref{prop: cyclically oriented lemma}, its weight is $1$ as well. Thus all weights of $Q$ are $1$ in either case, and \eqref{eq:unit-weight-classification} applies.
	\end{proof}

\medskip\noindent{\bf Acknowledgements}.
The first author is supported by JSPS KAKENHI Grant Number 25KJ1438.
The second author was supported by JSPS KAKENHI Grant Number JP21K13761.

\bibliographystyle{alpha}
\bibliography{bib}

\end{document}